\documentclass[10pt]{amsart} \usepackage{amssymb, verbatim, epic, eepic}

\def\Im{{\hbox{Im}}}
\def\sgn{{\hbox{sgn}}}
\def\Re{{\hbox{Re}}}
 
\def\dist{{\hbox{\rm dist}}}
\def\C{{\hbox{\bf C}}}
\def\R{{{\mathbf{R}}}} 

\def\Hess{{\operatorname{Hess}}}
\def\eps{\varepsilon}

\def\emph#1{{\it #1}} \def\textbf#1{{\bf #1}}
 
\newcommand{\nabb}{\mbox{$\nabla \mkern-13mu /$\,}}

\def\RR{\mathbb{R}}

\def\near{{\operatorname{near}}} 
\def\asympt{{\operatorname{asympt}}}
\def\Anear{A_{\near}} 
\def\Aasympt{A_{\asympt}} 
\def\ang{{\operatorname{ang}}}
\def\para{\operatorname{par}}
\def\Tscstar{{}^{\sc} T^*} 
\def\sc{\operatorname{sc}}
\def\Psisc{\Psi_{\sc}} 
\def\Id{\operatorname{Id}}

\newcommand{\pa}{\partial}
\newcommand{\Mbar}{\overline{M}}
\newcommand{\abs}[1]{{\left\lvert{#1}\right\rvert}}
\newcommand{\norm}[1]{{\left\lVert{#1}\right\rVert}}
\newcommand{\angs}[1]{{\left\langle{#1}\right\rangle}}
\newcommand{\ep}{{\epsilon}}
\newcommand{\CI}{{\mathcal{C}^\infty}}

\newcommand\lb{\operatorname{lb}}
\newcommand\rb{\operatorname{rb}}
\newcommand\blf{\operatorname{bf}}
\newcommand\conic{{\operatorname{conic}}}
\newcommand\CIdot{\dot C^\infty}
\DeclareMathOperator{\Tr}{{Tr}}
\DeclareMathOperator{\Op}{{Op}}
\DeclareMathOperator{\supp}{{supp}}

\theoremstyle{plain} \newtheorem{theorem}[subsection]{Theorem}
  
  \newtheorem{proposition}[subsection]{Proposition}
  \newtheorem{lemma}[subsection]{Lemma}
  \newtheorem{corollary}[subsection]{Corollary}

\theoremstyle{remark} 
  \newtheorem{remark}[subsection]{Remark}
  
  \newtheorem{example}[subsection]{Example}

\theoremstyle{definition} \newtheorem{definition}[subsection]{Definition}
  
\numberwithin{equation}{section}

\begin{document}

\title[Strichartz on manifolds]{A Strichartz inequality for the
Schr\"odinger equation on non-trapping asymptotically conic manifolds}

\author{Andrew Hassell} \thanks{A.H.\ is supported in part by an Australian
Research Council Fellowship.}  \address{Department of Mathematics, ANU,
Canberra, ACT 0200, AUSTRALIA} \email{hassell@maths.anu.edu.au}

\author{Terence Tao} \thanks{T.T.\ is a Clay Prize Fellow and is supported
in part by grants from the Packard Foundation.}  \address{Department of
  Mathematics, UCLA, Los Angeles California 90095, USA} \email{tao@math.ucla.edu}

\author{Jared Wunsch} \thanks{J.W.\ is supported in part by NSF grant
DMS-0323021.}  \address{Department of Mathematics, Northwestern University,
  Evanston IL 60208, USA}
\email{jwunsch@math.northwestern.edu}

\thanks{The authors thank an anonymous referee for helpful comments on the manuscript.}



\vspace{-0.3in}
\begin{abstract}
We obtain the Strichartz inequality
$$ \int_0^1 \int_M |u(t,z)|^4 \ dg(z) dt \leq C \| u(0) \|_{H^{1/4}(M)}^4$$
for any smooth three-dimensional Riemannian manifold $(M,g)$ which is
asymptotically conic at infinity and non-trapping, where $u$ is a solution
to the Schr\"odinger equation $iu_t + \frac{1}{2} \Delta_M u = 0$.  The
exponent $H^{1/4}(M)$ is sharp, by scaling considerations.  In particular
our result covers asymptotically flat non-trapping manifolds.  Our argument
is based on the interaction Morawetz inequality introduced by Colliander et
al., interpreted here as a positive commutator inequality for the tensor
product $U(t,z',z'') := u(t,z') u(t,z'')$ of the solution with itself. We
also use smoothing estimates for Schr\"odinger solutions including 
one (proved here) with weight $r^{-1}$ at infinity and with the gradient term
involving only \emph{one} angular derivative.
\end{abstract}

\maketitle

\tableofcontents

\section{Introduction}\label{sec:introduction}

The purpose of this paper is to establish a $L^4$ space-time Strichartz
inequality on a class of non-Euclidean spaces, namely smooth
three-dimensional asymptotically conic Riemannian manifolds
$(M,g)$ which obey a non-trapping condition.  We make these concepts
precise as follows.

\begin{definition}\label{M-def}  A smooth complete noncompact Riemannian manifold $M$, parameterized\footnote{We use $z$ as the spatial co-ordinate on $M$, in order
to reserve the letter $x$ for the scattering co-ordinate $x := 1/r$.  While
we phrase some of these definitions in general dimension $n$, we shall
focus primarily on the three-dimensional case $n=3$.  We will use the
indices $j,k$ to parameterize three-dimensional manifolds $M$ or
two-dimensional manifolds $\partial M$ (reserving $i$ for the square root
of -1), and the indices $\alpha, \beta$ to parameterize the six-dimensional
product manifold $M \times M$ which will play a prominent role in the
argument.}  by the variable $z := (z^1, \dots, z^n)$, with metric
$g:=g_{jk}(z) dz^j dz^k$, is \emph{asymptotically conic} if there exists a
fixed compact set $K_0 \subset M$ and an $(n-1)$-dimensional smooth compact
Riemannian manifold $(\partial M, h)$ --- parameterized by the variable $y
:= (y^1, \dots, y^{n-1})$, with metric $h := h_{jk}(y) dy^j dy^k$ --- such
that the \emph{scattering region} or \emph{asymptotic region} $M \backslash
K_0$ can be parameterized as the collar neighbourhood
$$ M \backslash K_0 \equiv (0, \epsilon_0) \times \partial M := \{ (x, y):
0 < x < \epsilon_0, y \in \partial M \}$$ for some $\epsilon_0 > 0$, with
metric of the form
\begin{equation} g_{jk}(z) dz^j dz^k = \frac{dx^2}{x^4} + \frac{h_{jk}(x,y) dy^j
dy^k}{x^2}
\label{ac-metric}\end{equation}
 with the usual summation conventions; here $h_{jk}(x,y)$ is a smooth
function on $[0,\epsilon_0) \times \partial M$ such that $h_{jk}(0,y) =
h_{jk}(y)$. We will often use the coordinate $r := 1/x$ instead of $x$ in
the scattering region.\footnote{Note that $\pa M$ can in fact be recovered
  from the Riemannian manifold $M$ by identifying points in $\pa M$ with
  pencils of geodesics that remain at bounded distance apart as $t\to\infty$.}

The metric in
co-ordinates $(r, y)$ has the `asymptotically conic' form
$$ g = dr^2 + r^2 h_{jk}(\frac{1}{r}, y) dy^j dy^k.
$$ By analogy with Euclidean space, we call $r$ the \emph{radial variable},
and $y$ the \emph{angular variable}; the latter variable is often also
denoted $\omega$ or $\theta$ in other texts; here we shall reserve Greek
letters for `cotangent' or `frequency' co-ordinates and Latin letters
for `spatial' or `position' co-ordinates.  For $V \in T_z M,$ we let
$\abs{V}_{g(z)}$ denote its length with respect to the metric $g.$

In particular we may compactify $M$ to $\overline M := M \cup \partial M$
by identifying $\partial M$ with $\{0\} \times \partial M$ in this
co-ordinate chart.  If $h_{ij}$ is independent of $x$ (so that $h_{ij}(x,y)
= h_{ij}(y)$) we say that $M$ is \emph{perfectly conic near infinity}.
Finally, we say that $M$ is \emph{non-trapping} if 
every geodesic $z(s)$ in $M$ reaches $\partial M$ as $s \to \pm \infty$.

It will be convenient to select a somewhat artificial, but globally defined radial positive weight $\langle z \rangle$,
chosen so that $\langle z \rangle$ is equal to $r$ in the scattering region $(0,\epsilon_0) \times \partial M$
and is comparable to $1/\epsilon_0$ in the compact interior region $K_0$, in such a way that
the function $\langle z \rangle$ is smooth and obeys the symbol estimates
$$ |\nabla_z^j \langle z \rangle| \leq C_j \langle z \rangle^{1-j}$$
for all $j \geq 0$.
\end{definition}

\begin{remark}
The most important example of an asympotically conic manifold is Euclidean
space $(M,g) := (\R^n,\delta)$, with $K_0$ equal to the unit ball $K_0 := \{ z:
|z| \leq 1 \}$, say, with $\partial M$ equal to the unit sphere
$S^{n-1} = \{ z: |z| = 1 \}$ with its standard metric $h$,
and scattering co-ordinates
$$ x := \frac{1}{|z|}; \quad y := \frac{z}{|z|}; \quad r := |z|$$ for $r >
r_0 := 1$.  This example is in fact perfectly conic near infinity, and is
also clearly non-trapping.  More generally, any compact perturbation of
Euclidean space will be perfectly conic near infinity, though it may not be
non-trapping.  Any asympotically
Euclidean space $(\R^n, g)$ with decay estimates $|\nabla^j (g-\delta)(z)|
\leq C_j \langle r\rangle^{-j-1}$ will also be asymptotically conic, with
$\partial M$ equal to the standard sphere $S^{n-1}$; it is quite likely
that the smoothness (and decay) assumptions we use can be weakened
substantially.

\begin{remark} Note that a non-trapping manifold must in fact be contractible.  To see
this, we observe that the region $K:= \{x\geq \ep\}$ is, for $\ep$
sufficiently small, a manifold with the same homotopy type as $M$ which is
\emph{convex} in the sense required to apply Theorem~4.2 of
\cite{Thorbergsson}.  This guarantees the existence of a closed geodesic in
$M$ provided $\pi_j (K)\neq 0$ for some $j>0$. \end{remark}

\end{remark}

\begin{remark}
We note that the class of metrics defined above is precisely the class of
\emph{scattering metrics} defined by Melrose \cite{melrose}, with the
metric written in a normal form in the scattering region due to Joshi and
S\'a-Barreto \cite{jsb}.
\end{remark}

Let $(M, g)$ be an asymptotically conic manifold, and suppose $u$ is a smooth
solution to the (time-dependent) Schr\"odinger equation
\begin{equation}\label{nls}
iu_t + \frac{1}{2} \Delta_M u = 0
\end{equation}
on $\R \times M$, where $u(t,z)$ is a complex field, $u(0) \in \CIdot(\Mbar)$ is a Schwartz function (i.e.\ it and its
derivatives vanishes to infinite order at $\partial M$)
 and
\begin{equation}\label{beltrami}
\Delta_M := \frac{1}{\sqrt{g}} \partial_j \sqrt{g} g^{j k}
\partial_k = \Tr_g \nabla^2, \quad g := \det(g_{jk})
\end{equation}
is the (negative definite) Laplace-Beltrami operator (with $\nabla$
denoting covariant derivatives).  It is well known
(see e.g.\ \cite{cks}) that if $u(0)$ is Schwartz then there is a unique global Schwartz solution to \eqref{nls},
but our estimates will not depend on any Schwartz norms of $u$.  We may rewrite \eqref{nls} as
\begin{equation}\label{nls-H}
 u_t = -iH u,
\end{equation}
where $H := - \frac{1}{2} \Delta$ is the Hamiltonian; note that
$H$ is positive definite and self-adjoint with respect to the real-valued
inner product
\begin{equation}\label{inner-def}
 \langle u, v \rangle_M := \Re \int_M u(z) \overline{v(z)}\ dg(z)
\end{equation}
where $dg(z) := \sqrt{g} \, dz$ is the usual measure induced by the Riemannian metric $g$.  In particular
we can generate the usual functional calculus of $H$, and define Sobolev spaces $H^s(M)$ for all $s \in \R$
as the class of functions $u$ whose $H^s(M)$ norm $\| u \|_{H^s(M)} := \| (1+H)^{s/2} u \|_{L^2(M)}$ is 
finite, where of course we use \eqref{inner-def} to define $L^2(M)$.  

The purpose of this paper is to establish the following local-in-time $L^4$
Strichartz estimate:

\begin{theorem}\label{main}  For any smooth three-dimensional asymptotically conic non-trapping manifold $M$,
and any (Schwartz) solution to \eqref{nls} (or \eqref{nls-H}), we have
\begin{equation}\label{l4-strichartz}
 \int_0^1 \int_M |u(t,z)|^4\ dg(z) dt \leq C \| u(0) \|_{H^{1/4}(M)}^4.
\end{equation}
\end{theorem}

Here and in the sequel constants such as $C$ are allowed to depend on the Riemannian 
manifold $(M, g)$, but not on $u$.  In particular once one proves this estimate
for Schwartz initial data $u(0)$, one can obtain the same estimate for
general $H^{1/4}(M)$ data by the usual limiting argument.

We now contrast Theorem~\ref{main} with existing results.  In Euclidean
space $(M, g) = (\R^n, \delta)$ we have the \emph{Strichartz estimate} 
(see e.g.\ \cite{tao:keel} and the references therein)
\begin{equation}\label{strich}
\| u \|_{L^q_t L^r_x(\R \times \R^n)} \leq C_{n,q,r,s} \| u(0) \|_{H^s(\R^n)}
\end{equation}
whenever
\begin{equation}\label{cond}
 \frac{n}{2} - s \leq \frac{2}{q} + \frac{n}{r} \leq \frac{n}{2}; \quad 2
 \leq q \leq \infty; \quad 2 \leq r < \infty.
\end{equation}
Here of course
$$
\| u \|_{L^q_t L^r_x(\R \times \R^n)} := \left(\int_\R \left(\int_{\R^n}
    |u(t,z)|^r\ dz\right)^{q/r}\ dt\right)^{1/q}.$$
In particular 
the estimate \eqref{l4-strichartz} is the case $(n,q,r,s) =
(3,4,4,1/4)$, but only locally in time: $t \in [0,1]$. 
In Euclidean space a scaling argument, replacing $u(t,z)$ by $u(\lambda^2 t, \lambda z)$ and
letting $\lambda \to \infty$,  shows that the
estimate \eqref{l4-strichartz} holds globally in time, that the exponent $H^{1/4}$ on the
right-hand side of \eqref{l4-strichartz} is sharp, and that we can even replace
the inhomogeneous Sobolev space $H^{1/4}(\R^3)$ by the homogeneous $\dot
H^{1/4}(\R^3)$.  We believe that the global-in-time estimate also holds in
the asymptotically conic nontrapping setting, but we have not attempted to prove it here.   

Strichartz estimates are crucial in analyzing the low regularity behavior of non-linear
Schr\"odinger equations, see for instance \cite{borg:book}, \cite{caz}, and
the remarks in Section \ref{remarks-sec}.

In Euclidean space there is also a \emph{local smoothing estimate}
\begin{equation}\label{local-smoothing}
 \int_\R \int_K |\nabla u(t,z)|^2\ dx dt \leq C_K \| u(0)
 \|_{H^{1/2}(\R^3)}^2
\end{equation}
for any compact set $K \subset \R^3;$ such estimates seem to have first
appeared in the work of Constantin and Saut \cite{CS}, Sj\"olin
\cite{sjolin} and Vega \cite{vega}.  (Earlier estimates of this
`dispersive smoothing' type, for the KdV equation, date back to work of
Kato \cite{Kato}.)  In Euclidean space the estimates \eqref{strich},
\eqref{local-smoothing} are not directly related---indeed they can be
proven independently of each other.  However, on other manifolds the local
smoothing estimate seems to be indispensable for proving optimal Strichartz
estimates.

For general manifolds $(M,g)$, the local smoothing estimate
is in general false, even if we localize in time, if there are trapped geodesics; this was shown by Doi \cite{doi}.  However,
assuming that $M$ is non-trapping, we have the analogue of
\eqref{local-smoothing} in arbitrary dimension:
\begin{equation}\label{local-smoothing-manifold}
 \int_0^1 \int_K |\nabla u(t,z)|_{g(z)}^2\ dg(z) dt \leq C_K \| u(0)
 \|_{H^{1/2}(M)}^2.
\end{equation}
This estimate is due to Craig, Kappeler and Strauss \cite{cks} on
asymptotically Euclidean space (see also \cite{doi}).  We shall need this estimate as well as several global versions of it. Some well-known versions (at least for asymptotically flat manifolds; the generalization to asymptotically conic is straightforward) are 
\begin{equation}
\label{eq:manifold-doi-comp}
\int_0^1 \int_M \frac{|\nabla u(t,z)|_{g(z)}^2}{\langle z \rangle^{1+\eps}}\ dg \, dt \leq C_\eps \| u(0)
\|_{H^{1/2}(M)}^2
\end{equation}
as well as the variant
\begin{equation}
\label{eq:manifold-doi-noderivs}
\int_0^1 \int_M \frac{|u(t,z)|^2}{\langle z \rangle^{1+\eps}}\ dg \, dt \leq C_\eps \norm{u(0)}_{H^{-1/2}(M)}^2
\end{equation}
for $\eps > 0$.  These are proved in Appendix II. Furthermore, if $r_0 > 0$ is sufficiently large, we have the Morawetz estimate
\begin{equation}
\label{eq:manifold-doi-ang}
\int_0^1 \int_{\langle z \rangle > r_0} \frac{|\nabb u(t,z)|_{g(z)}^2}{\langle z \rangle}\ dg(z) \, dt \leq C \| u(0)
\|_{H^{1/2}(M)}^2
\end{equation}
where $\nabb := x \nabla_y= \angs{z}^{-1} \nabla_y,$ which is well-defined
in the scattering region, and thus on the domain of integration if $r_0$ is
sufficiently large.  Note that for the purposes of this estimate, we could
equally well define $\nabb$ as the orthogonal projection of $\nabla$ to the
hyperplane in $T_z M$ orthogonal to $\pa_r,$ since $x\nabla_{y_i}$ form a
basis for this hyperplane.  (For simplicity in writing our estimates, we
will henceforth take $\nabb$ to be defined globally, and equal to $0$ in
the compact region $\angs{z}<r_0.$)

In this paper, we prove the following generalization of \eqref{eq:manifold-doi-ang}:

\begin{lemma}\label{superDoi-intro}
Let $M$ be an asymptotically conic non-trapping manifold. Let $a^{jk}(z)$
be symbols of order zero on $M$, i.e. functions obeying the estimates 
$$ |\nabla_z^m a^{jk}(z)| \leq C_{m} \langle z \rangle^{-m}$$
for all $m \geq 0$. If $u(0) \in H^{1/2}(M)$ and $u$ is the corresponding solution to \eqref{nls}, then 
\begin{equation*}
\int_0^1  \bigg| \int_{M} \frac{a^{jk}(z) \nabb_j u(t,z)
\overline{\nabla_k u(t,z) }}{r} \, dg(z) \bigg| dt \leq C \| u(0)
\|_{H^{1/2}(M)}^2.
\end{equation*} 
\end{lemma}

This improves upon \eqref{eq:manifold-doi-ang} since only one of the $u$-derivatives is required to be angular. In flat Euclidean space it was proved recently by Sugimoto \cite{Sug}. 
In fact, we prove (and need) a stronger version of this result where the symbols $a^{jk}$ are allowed to depend on an additional parameter; see Lemma~\ref{lemma:superDoi}.



The epsilon loss in \eqref{eq:manifold-doi-comp},
\eqref{eq:manifold-doi-noderivs} is not removable, even in $\RR^n$
(although refinements with logarithmic losses are certainly possible);
semi-classically, this can be explained by the fact that a particle $z(s)$
moving at unit speed along a geodesic will have a finite integral $\int_\R
\langle z(s) \rangle^{-1-\eps}\, ds$, whereas the integral $\int_\R \langle
z(s) \rangle^{-1}\, ds$ diverges logarithmically.  The Morawetz estimate
\eqref{eq:manifold-doi-ang} then corresponds semi-classically to the fact
that the slightly smaller integral $\int_{\langle z(s) \rangle > r_0}
|v_\ang(s)|^2 \langle z(s) \rangle^{-1}\, ds$ avoids the logarithmic
divergence and converges absolutely, where $v_\ang(s) = r \frac{d}{ds}
y(s)$ is the angular component of the velocity.  This reflects the fact
that the velocity vector of a geodesic becomes purely radial in the
asymptotic limit.

Now consider the Strichartz estimate \eqref{strich} on general manifolds.
Without the non-trapping condition, one does not have any local smoothing,
and one does not expect to obtain the estimate \eqref{strich} with the
sharp number of derivatives $s$.  However, it is still possible to obtain a
(local-in-time) Strichartz estimate with a loss of derivatives.  Indeed,
Burq, Gerard and Tzvetkov \cite{bgt} showed that the estimate
\begin{equation}\label{strich-manifold}
\| u \|_{L^q_t L^r_x([0,1] \times M)} \leq C \| u(0)
 \|_{H^{s + \frac{1}{q}}(M)}
\end{equation}
holds for arbitrary smooth manifolds (trapping or non-trapping) in general dimension.
 This loss of derivatives, as well as the localization in time, 
was shown in \cite{bgt} to be sharp in the case of the sphere (which is in some sense maximally trapping).

The next progress was by Staffilani and Tataru \cite{st}, who were able to
remove the $1/q$ loss of derivatives for metrics on $\R^n$ which are non-trapping and Euclidean
outside a compact set (and also only required $C^2$ regularity of the metric); a key tool was the use of
\eqref{local-smoothing-manifold} to localize to a compact region of space.
After the work of Staffilani and Tataru, Burq \cite{burq} gave an
alternative proof of the same result, and also considered metrics which were asymptotically flat rather than flat outside of a compact set.  In this more general setting Burq was almost able to recover the
same Strichartz estimate as Staffilani and Tataru, but with an epsilon
loss\footnote{The endpoint $q=2$ is also not obtained, for a rather
technical reason involving the causality cutoff $s < t$ in Duhamel's
formula.} of derivatives:
$$ \| u \|_{L^q_t L^r_x([0,1] \times M)} \leq C \| u(0)
 \|_{H^{s + \eps}(M)}.
 $$ 
 The idea was to divide $M$ into dyadic shells $R
 \leq |z| \leq 2R$ and apply a variant of the Staffilani-Tataru argument on
 each shell, either rescaling $R$ to become 1, or relying on
\eqref{eq:manifold-doi-comp} instead of \eqref{local-smoothing-manifold}.
As remarked earlier, the presence of
the $\eps$ in \eqref{eq:manifold-doi-comp} is
essential, and is related to the epsilon
loss of derivatives in Burq's result.

Our result in Theorem \ref{main} is slightly stronger than Burq's in the
sense that it applies also to asymptotically conic manifolds, and removes
the epsilon loss.  However, it is also weaker than these previous results
because it is restricted to the specific exponents $(n,q,r,s) =
(3,4,4,1/4)$ and cannot obtain the full range \eqref{cond} of the estimate
\eqref{strich-manifold}, although one can obtain a subset of this range by
applying Sobolev embedding to \eqref{l4-strichartz} and also using the
energy identity $\| u(t) \|_{L^\infty_t L^2(\R \times M)} = \| u(0)
\|_{L^2(M)}$.  We conjecture however that the full set of Strichartz
estimates can be recovered with no loss of derivatives for non-trapping
asymptotically conic (or asymptotically flat) manifolds in arbitrary
dimension.

Unlike the previous results, the proof of Theorem \ref{main} does not
proceed via construction of a parametrix or proving any dispersive
estimates on the fundamental solution.  Instead, we use the interaction
Morawetz inequality approach introduced in \cite{ckstt:scatter},
\cite{ckstt:scatter-2} in the context of a non-linear Schr\"odinger
equation in Euclidean space.  We re-interpret this approach in the language
of positive commutators applied to the tensor product $U(t,z',z'') :=
u(t,z') u(t,z'')$ of the solution $u$ with itself, and then modify it to
asymptotically conic manifolds.

There will be a number of technical difficulties in carrying out the above
strategy, first in understanding the error terms generated by the non-zero
curvature, and second the problem of avoiding the singularities of the
metric function $d_M(z',z'')$ once this distance exceeds the radius of
injectivity.  Unlike the argument in \cite{burq}, our approach cannot afford to
localize the radial variable to dyadic blocks $R < r < 2R$ as this would
create error terms with radial derivatives in them, which cannot be
controlled by \eqref{eq:manifold-doi-comp} without losing an epsilon too many
derivatives.  In fact we can only localize $r$ once, to the region $r >
r_0$, though we can localize the \emph{angular} variable $y$ to a small
set
\footnote{Even doing this localization creates some unpleasant error terms in which only \emph{one} of the derivatives is angular. We shall control such terms using Lemma~\ref{superDoi-intro}, or more precisely Lemma~\ref{lemma:superDoi}, which only requires one angular derivative.}.  
The understanding of the
geometry in this regime --- in particular, controlling the derivatives of the
distance function $d_M((r',y'), (r'',y''))$ when $y', y''$ are close, and
$r', r''$ are both large but not necessarily comparable --- appears to be
absolutely essential for obtaining bounds such as \eqref{strich} with no
loss of derivatives whatsoever.  Thus our argument requires a certain
amount of geometric machinery which we have placed in an appendix
(Section \ref{appendix}).

It is well known that Strichartz estimates can be used in the theory of
nonlinear Schr\"odinger equations, and we briefly discuss a modest application of this
estimate to a class of non-linear problems in Section \ref{remarks-sec}.  
However as our Strichartz estimate is restricted to $L^4$ in space and time, the application is admittedly somewhat limited.  It is still an important
problem as to whether one can obtain the full range of Strichartz estimates (with no derivative
loss) on general non-trapping manifolds, and whether they can be extended to be global in time; 
such a result will have substantially more applications
to non-linear Schr\"odinger equations.



\section{Positive commutators and Morawetz identities}\label{morawetz-sec}

In this section we review the basic method of positive commutators in the
context of the time-dependent Schr\"odinger equation, and also indicate how
it is connected to estimates of Morawetz type; we will rely on these
methods, together with detailed analysis of the various terms arising from
such methods, to prove Theorem \ref{main}.  In this discussion $M$ is a
general smooth Riemannian manifold with compactification $\Mbar$.  In our
applications $M$ will either be a three-dimensional non-trapping
asymptotically conic manifold, or else the product manifold $M \times M$
which is six-dimensional and non-trapping, but not asymptotically conic (as
we shall see, it has a `corner' at infinity which we will have to blow up
to analyze).  To emphasize that we may be working in the product setting,
we will use Greek indices $\alpha,\beta,\gamma$ to parameterize the
manifold rather than Latin indices $j,k,l$.

Let $u$ be a solution to the Schr\"odinger equation \eqref{nls} (or
\eqref{nls-H}); to avoid technical issues let us assume that $u \in \CIdot(\Mbar)$, that is, $u$ is rapidly decreasing together with all derivatives as $z \to \infty$.   Given a pseudo-differential\footnote{In fact, we shall mostly be able to rely on \emph{classical} differential operators $A$ for our commutants; pseudo-differential commutants will only be needed to prove the local smoothing estimates necessary to control errors arising from the classical
commutant.} operator $A$ on $M$, we can easily verify\footnote{Note the presence of the real part in \eqref{inner-def} allows us to obtain this even when $A$ is not self-adjoint.  In practice, we will only deal with operators $A$ whose \emph{principal} symbols are real, so that they are self-adjoint modulo lower order terms.} the \emph{Heisenberg
equation}
$$ \partial_t \langle A u(t), u(t) \rangle_M = \langle i[H, A]
u(t), u(t) \rangle_M,
$$ 
where $[H,A] := HA-AH$ is the usual commutator.  Integrating this we obtain the inequality
\begin{equation}\label{heisenberg}
\int_0^1 \langle i[H, A] u(t), u(t)\rangle_M\ dt \leq 2 \sup_{t \in [0,1]}
|\langle Au(t), u(t) \rangle_M|.
\end{equation}
The strategy of the method of positive commutators is to select a good commutant
$A$, chosen so that the commutator $i[H, A]$ is positive
(semi-)definite, possibly modulo a lower order error term, and then apply
\eqref{heisenberg} to obtain a non-trivial spacetime bound on $u$.  From 
the pseudo-differential calculus and the sharp G\aa rding inequality, one only
expects to be able to do this if the symbol $\sigma(A)$ of $A$ increases
along the bicharacteristic flow of $H$ (i.e.\ geodesic flow in phase space).
At first glance this method only seems able to obtain quadratic estimates
on $u$ (e.g.\ the method can be used to deduce \eqref{local-smoothing-manifold} -- \eqref{eq:manifold-doi-ang}), but we will also be able to obtain quartic estimates (and in
particular the $L^4$ bound \eqref{l4-strichartz}) by the trick of replacing $u$ by the
tensor product $U(t,z',z'') := u(t,z') u(t,z'')$, and working in the 
six-dimensional product manifold $M \times M$.

For future reference we compute some commutators of $H$ with multiplier operators.  

\begin{lemma}\label{commute}  Let $a(x)$ be a real-valued tempered distribution on $M$, thought of as a multiplier operator $(af)(x) := a(x) f(x)$ on Schwartz functions.  Then we have the single commutator identity
\begin{equation}\label{a-bracket}
i[H, a] = i \langle (\nabla a), \nabla \rangle_g + i (H a) 
= i (\nabla^\alpha a) \nabla_\alpha + i (Ha)
\end{equation}
and the double commutator identity
\begin{equation}\label{double-bracket}
-[H, [H, a]] = - \nabla_\beta \operatorname{Hess}(a)^{\alpha \beta} \nabla_\alpha - (H^2 a)
\end{equation}
where $\operatorname{Hess}(a)^{\alpha\beta}$ is the symmetric tensor
\begin{equation}\label{Q-def}
\operatorname{Hess}(a)^{\alpha \beta} = (\nabla da)^{\alpha \beta} =
g^{\alpha \gamma} g^{\beta \delta} \Big( \partial_\gamma \partial_\delta a + \Gamma^\rho_{\gamma \delta} \partial_\rho a \Big).
\end{equation}
\end{lemma}

\begin{proof}
The proof of \eqref{a-bracket} is elementary. To prove \eqref{double-bracket}, observe that both sides of
\eqref{double-bracket} are second-order self-adjoint operators (with
respect to $\langle, \rangle_M$) with real coefficients.  Also, by applying
both sides to the constant function 1 we see that the constant terms of
\eqref{double-bracket} match.  Thus it will suffice to show that the
principal symbols of both sides match, since this would mean that the
difference of the left-hand side and right-hand side is a self-adjoint
first-order operator with real coefficients and no constant term, which is
necessarily zero.

To compute the principal symbols at a point $p$, we pass to Riemannian
normal coordinates about $p$.  Then up to second order at $p$, $g_{\alpha\beta}=\delta_{\alpha\beta},$ and inserting
\eqref{a-bracket}, we see that the LHS of \eqref{double-bracket} becomes
$$
-(\pa_\alpha\pa_\beta(a))\pa_\alpha\pa_\beta,
$$ modulo lower order terms.  That the symbol agrees with minus the Hessian
follows, as Christoffel symbols in \eqref{Q-def} vanish at $p.$

\end{proof}

\begin{remark}  Because $i[H,a]$ is self-adjoint, we know that $\int
  \angs{i[H,a]u, u} \, dg(z)$ is real for any test function $u$.  Since
  $\int \angs{i(Ha)u, u}\, dg(z)$ is clearly imaginary, we thus see from \eqref{a-bracket} that
\begin{equation}\label{ad} \langle i[H,a] u, u\rangle_M = \Im \int_M \langle \nabla a, \nabla u \rangle_g \overline{u}\ dg. 
\end{equation}
More generally, if $b$ is another real function, then
\begin{equation}\label{adb} \langle b i[H,a] u, u\rangle_M = 
\Im \int_M b \langle \nabla a, \nabla u \rangle_g \overline{u}\ dg. 
\end{equation}
\end{remark}

\begin{remark}
The Hessian $\Hess(a)^{\alpha \beta}$ measures the convexity of the function
$a(z)$ with respect to the geodesic flow.  Indeed, if $z(s)$ is the curve
of a unit-speed geodesic in $M$, and $v^\alpha(s) := \frac{d}{ds}
z^\alpha(s)$ is the velocity vector, then one has
\begin{equation}\label{classical}
\frac{d^2}{ds^2} a(z(s)) = \Hess(a)_{\alpha \beta}(z(s)) v^\alpha(s)
v^\beta(s).
\end{equation}
This can be seen either from \eqref{Q-def} and the geodesic flow equation,
or alternatively from \eqref{double-bracket} and the symbol calculus.
\end{remark}

\begin{example}[Euclidean one-particle Morawetz
    inequality; \cite{moraw}, \cite{ls}]\label{1-particle}  Consider flat
three-space $M := \R^3$.  Let $a(z) := |z|$, the distance function to the
origin. A simple computation shows that
$$ \Hess(a)^{jk}(z) = \partial_j \partial_k a(z) =
\frac{1}{|z|} - \frac{z^j z^k}{|z|^3}, \text{ so that}
$$
$$ \Hess(a)^{jk}(z) \overline{\partial^j u(t,z)} \partial^k
u(t,z) = \frac{|\nabb u(t,z)|^2}{|z|}$$ 
where $\nabb$ denotes the angular
gradient as before.  Also, 
$$ -H^2 a = -\frac{1}{4} \Delta\Delta a = 2 \pi \delta,$$ where $\delta$ is
the Dirac mass at the origin in $\R^3$. Thus by
\eqref{double-bracket}, \eqref{ad} we have
$$ \langle -[H, [H, a]] u, u \rangle_M = 2 \pi |u(t,0)|^2 + \int_{\R^3}
\frac{|\nabb u(t,z)|^2}{|z|}\ dz$$ and
$$ \langle i[H, a] u, u \rangle_M = \Im \int \overline u(t,z) \frac{z}{|z|}
\cdot \nabla u(t,z)\ dz.$$ 
Since $\nabla$ maps $\dot H^{\frac1{2}}(\R^3)$ to $\dot H^{-\frac1{2}}(\R^3)$, and
$\frac{z}{|z|}$ is a bounded multiplier\footnote{This can be seen by first 
verifying boundedness on $L^2(\R^3)$ and $\dot H^1(\R^3)$ (using Hardy's
inequality) and then using interpolation and duality.} on $\dot H^{-\frac1{2}}(\R^3)$
we have
$$ \abs{\int_{\R^3} \overline u(t,z) \frac{z}{|z|} \cdot \nabla u(t,z)\ dz} \leq C
\| u(t) \|_{\dot H^{1/2}(\R^3)}^2;$$
since the $\dot H^{1/2}(\R^3)$ norm is
conserved, setting $A=i[H,a]$ in \eqref{heisenberg} we thus obtain the
\emph{Morawetz inequality}
$$ \int_0^1 |u(t,0)|^2\ dt + \int_0^1 \int_{\R^3} \frac{|\nabb
u(t,z)|^2}{|z|}\ dz dt \leq C \| u(0) \|_{\dot H^{1/2}(\R^3)}^2.$$ 
In fact,
in Euclidean space this particular estimate extends from the time interval
$[0,1]$ to all of $\R$ by the same argument.  Note in particular that we have proven \eqref{eq:manifold-doi-ang}
for Euclidean spaces.  A similar argument can be used to derive \eqref{eq:manifold-doi-ang} for general
non-trapping manifolds, as well as \eqref{eq:manifold-doi-comp}, if one
already knows the dispersive smoothing estimate \eqref{local-smoothing-manifold};
however the estimate \eqref{local-smoothing-manifold} cannot in general be proven without introducing commutants
which are pseudo-differential operators rather than classical differential operators.
\end{example}

\begin{example}[Euclidean two-particle (interaction) Morawetz inequality; \cite{ckstt:scatter}, \cite{ckstt:scatter-2}]\label{2-particle}
We now consider a variant of the Morawetz inequality.  We let $M$ be the
six-dimensional space $M := \R^3 \times \R^3$, which we parameterize as $z
= (z',z'')$ where $z',z'' \in \R^3$ and $z \in \R^3 \times \R^3$; observe
that the Hamiltonian $H := -\frac{1}{2} \Delta_{\R^3 \times \R^3}$ on this
manifold splits as $H = H_{z'} + H_{z''}$ where $H_{z'} := -\frac{1}{2}
\Delta_{z'}$ and $H_{z''} := -\frac{1}{2} \Delta_{z''}$.

Let $a: \R^3 \times \R^3 \to \R$ be the Euclidean distance function $a(z',z'') = |z'-z''|$. 
The function $a(z',z'')$ is clearly (weakly) convex under geodesic flow of
$z'$ and $z''$ (as can be seen for instance using center of mass co-ordinates), thus
the tensor $\Hess(a)^{\alpha \beta}$ is positive semi-definite.  In particular
we have
\begin{equation}\label{psd}
 \Hess(a)^{\alpha \beta}(z) \overline{\partial^\alpha U(t,z',z'')}
\partial^\beta U(t,z',z'') \geq 0.
\end{equation}
 One can also compute
\begin{equation}\label{flat-delta-delta}
-H^2 a = 8 \pi \delta(z'-z''),
\end{equation}
thus $-H^2 a$ is a Dirac mass on the diagonal $\{ (z',z'') \in \R^3:
z' = z'' \}$ of $\R^3 \times \R^3$.  The same reasoning as in the previous example yields the product Morawetz inequality
$$ \int_0^1 \int_{\R^3} |U(t,z,z)|^2\ dz dt \leq C \| U(0) \|_{\dot
H^{1/2}(\R^3 \times \R^3)}^2.$$ 
This inequality holds for all solutions
$U(t,z',z'')$ to the Schr\"odinger equation on $\R^3 \times \R^3$.  If we
now specialize to the tensor product $U(t,z',z'') := u(t,z') u(t,z'')$ of
an $\R^3$ solution $u(t,z)$ of the Schr\"odinger equation\footnote{Note
here we are using a very special property of the Schr\"odinger with itself,
namely that tensor products of solutions to the Schr\"odinger equation
remain solutions to a Schr\"odinger equation (but now on the product
manifold of twice the dimension).  This property is not shared by other
equations such as the wave or Klein-Gordon equation, and it remains a
challenge to extend this type of argument to those equations.}, we obtain
\begin{equation}\label{interaction}
\int_0^1 \int_{\R^3} |u(t,z)|^4\ dz dt \leq C \| u(0) \|_{L^2(\R^3)}^2 \|
u(0) \|_{\dot H^{1/2}(\R^3)}^2;
\end{equation}
after a Littlewood-Paley decomposition (see Proposition \ref{LP} below), this will imply in
particular the Strichartz estimate \eqref{l4-strichartz}. Again in the
Euclidean space setting it is possible to extend the time interval from
$[0,1]$ to all of $\R$.
\end{example}

We have just proven Theorem \ref{main}, for Euclidean space $\R^3$.  It is
thus tempting to adapt the above proof to more general manifolds.  The most
obvious thing to try is to repeat the above argument, with the manifold
$\R^3 \times \R^3$ replaced of course by $M \times M$ (with the product
metric $g \oplus g$), and with the function $|z'-z''|$ replaced by the
distance function $d_M(z',z'')$.  This works well locally, when $z'$ and
$z''$ are close together; the main issue is to verify that $d_M(z',z'')$
remains essentially geodesically convex, and that the error terms can be
treated by the one-particle estimates \eqref{eq:manifold-doi-comp} -- \eqref{eq:manifold-doi-ang}.  Difficulties arise however in the scattering region, requiring
more careful control on the geometry of the distance function on
asymptotically conic manifolds over long distances, as well as requiring a
refinement of Lemma~\ref{superDoi-intro}, but we shall address these
issues later.

In Euclidean space one can also develop $k$-particle versions of the
Morawetz inequality for any $k\geq 2$ by working in center-of-mass
coordinates, with the function $a$ being the standard deviation of the $k$
positions $z_1,\dots,z_k.$  This leads to a useful monotonicity formula in
$d$ dimensions when $d(k-1)\geq 3$ (for instance the $4$-particle Morawetz
inequality gives an $L^8_{t,z}$ estimate in one dimension).  In this paper,
however, we will confine ourselves to the two-particle problem.


\section{Littlewood-Paley reduction}\label{sec:morawetz}

In the previous section we discussed the method of positive commutators, and showed how it
can be used to prove $L^4_{t,x}$ type estimates such as \eqref{l4-strichartz}.  It is then
tempting to apply the same argument immediately to general manifolds, using the distance function 
$d_M(z',z'')$ as a substitute for $|z'-z''|$ as mentioned above.  Before we do so, however,
it is convenient to make a Littlewood-Paley reduction which allows us to restrict the solution $u$
to a single dyadic frequency range $[2^{2k-2}, 2^{2k+2}]$ in the spectrum of $H$.  This reduction is
standard (see e.g.\ \cite{sogge:wave}, Lemma 5.1 and the ensuing discussion).

\begin{proposition}\label{LP} Suppose we have
\begin{equation}
\int_0^1 \int_M |u(z,t)|^4 \, dg(z) \, dt \leq C \, 2^k \,  \| u(0) \|_{L^2(M)}^4,
\label{LP-est}\end{equation}
for all solutions $u(z, t)$ of \eqref{nls-H}  with $u(z,0)$ in the range of the spectral projection $\chi_{[2^{2k-2}, 2^{2k+2}]}(H)$, and all $k \geq 1$.  Then Theorem \ref{main} holds.
\end{proposition}

\begin{proof}
The solution to \eqref{nls-H} is given by $u(t) =
e^{-itH} u(0)$.  Furthermore we have a Littlewood-Paley decomposition
$$ \Id = \phi_0(H) + \sum_{k=1}^\infty \psi_k(H)$$ where $\phi_0$ is a
compactly supported bump function, and $\psi_k$ is a bump function adapted
to the interval $[2^{2k-2}, 2^{2k+2}]$.  This induces a decomposition
$$ u(0) = u(0)^{(0)} + \sum_{k=1}^\infty u(0)^{(k)}, \quad u(0)^{(0)} =
\phi_0(H) u(0), \quad u(0)^{(k)} = \psi_k(H) u(0). 
$$
Informally, $\psi_k(H)$
is the spectral projection to functions of frequency comparable to $2^k$.

The contribution of the low-frequency term $u(0)^{(0)}$ can be handled very easily using the Sobolev embedding $\| u \|_{L^4} \leq C \| \nabla u \|_{L^2}$ in four dimensions:
\begin{equation*}\begin{gathered}
\| e^{-itH} u(0)^{(0)} \|_{L^4_{t,x}([0,1] \times M)} \leq C \|
\nabla_{z,t} \big( e^{-itH} \phi_0(H) u(0) \big) \|_{L^2([0,1]\times M)}\\ \leq C \|
(1+H) e^{-itH} \phi_0(H) u(0) \|_{L^2([0,1] \times M)}  \leq C  \| u(0) \|_{L^2(M)} \leq \| u(0) \|_{H^{1/4}(M)}.
\end{gathered}\end{equation*}

So consider the high-frequency terms
\begin{equation}\label{hfreq} \| \sum_{k=1}^\infty e^{-itH} u(0)^{(k)} \|_{L^4_{t,x}([0,1] \times
M)}.
\end{equation}  We use the Littlewood-Paley estimate
$$ \| \sum_k \psi_k(H) f \|_{L^4(M)} \leq C \| (\sum_k |\psi_k(H)
f|^2)^{1/2} \|_{L^4(M)}$$ (see e.g.\ \cite{bgt} or \cite{Stein}, chapter VI, section 7.14; the point is that the heat
kernels $e^{-tH}$ or resolvents $(H-\lambda)^{-1}$ enjoy good decay and
regularity properties, and thus the Littlewood-Paley square function
operator is of Calder\'on-Zygmund type), to bound  \eqref{hfreq} by
$$ C \| (\sum_{k=1}^\infty |e^{-itH} u(0)^{(k)}|^2)^{1/2}
\|_{L^4_{t,x}([0,1] \times M)};$$ we can move the $l^2$ sum outside of the
$L^4$ norm to bound this by
$$ C \Big(\sum_{k=1}^\infty \| e^{-itH} u(0)^{(k)} \|_{L^4_{t,x}([0,1] \times
M)}^2 \Big)^{1/2}.$$ The estimate \eqref{LP-est} implies that \eqref{hfreq} is bounded by 
\begin{equation*} C \Big(\sum_{k=1}^\infty 2^{k/2} \| u(0)^{(k)} \|_{L^2(M)}^{2} \Big)^{1/2}.\end{equation*}
From the spectral localization of $u(0)^{(k)}$, we can estimate this by
\begin{equation*} C \Big(\sum_{k=1}^\infty \| (1 + H)^{1/8} u(0)^{(k)} \|_{L^2(M)}^{2} \Big)^{1/2}.\end{equation*}
Since $u(0)^{(k)} = \psi_k(H) u$, we see from spectral calculus that this expression is $O( \| (1+H)^{1/8} u(0) \|_{L^2(M)} ) = O( \| u(0) \|_{H^{1/4}} )$, as desired.
\end{proof}

Henceforth we fix $k > 0$ and write
$u = u^{(k)}$, thus $u$ is implicitly localized to the portion of the
spectrum of $H$ lying in the interval $[2^{2k-2}, 2^{2k+2}]$.  Our task is now to prove \eqref{LP-est}.
Observe from this spectral localization (and the unitarity of the operators $e^{-itH}$) that we have
the bounds
\begin{equation}\label{spectrum}
c 2^{ks} \| u(0) \|_{L^2(M)} \leq  \| u(t) \|_{H^s(M)} \leq C_s 2^{ks} \| u(0) \|_{L^2(M)}
\end{equation}
for all $s \in \R$ and $t \in [0,1]$.  In particular, we have
\begin{equation}\label{spectrum-2}
 \| \nabla u(t) \|_{L^2(M)} \leq C \| u(t) \|_{H^1(M)} \leq C 2^k \| u(0) \|_{L^2(M)}
\end{equation}
(the first inequality being easily verified by an integration by parts).  These inequalities give 
us a significant amount of freedom
to estimate various factors of $u$ in various Sobolev spaces, which will be very convenient in the argument
that follows.  Indeed, because of these estimates, the \emph{location} of the derivatives in our expressions to
be estimated will no longer be particularly relevant; what matters\footnote{Of course, weights such as $1/\langle z\rangle$ and $1/\langle z'\rangle$ will also be important to keep track of.} is the \emph{number} of such 
derivatives (which
range between 0 and 2, with the terms with two derivatives of course being the most difficult) as well as their
\emph{type} (angular derivatives will, in general, be easier to handle than radial derivatives, thanks for instance
to \eqref{eq:manifold-doi-ang}).

\section{Overview of proof of Theorem \ref{main}}\label{overview}

Now that we have localized $u$ to a single frequency range, we apply the
method of positive commutators discussed in Section \ref{morawetz-sec}.  We
will in fact use this method twice; once for the near region $K_0$, and
once the scattering region $M \backslash K_0$, using different commutants
for the two regions.

Let $0 < \eps \ll 1$ be a small constant to be chosen later.  We use $Z$ to denote the quantity
\begin{equation}\label{Z-def}
 Z := \eps^{-2} 2^k  \| u(0) \|_{L^2(M)}^4 + \eps \int_0^1 \int_{M} |u(z,t)|^4 \, dg(z) \, dt.
\end{equation}
Clearly, to prove \eqref{LP-est}, it suffices to show that
\begin{equation}\label{LP-Z} \int_0^1 \int_{M} |u(z,t)|^4 \, dg(z) \, dt = O(Z),
\end{equation}
since for $\eps$ small enough we can absorb the second term in \eqref{Z-def} into the left-hand side.

It remains to prove \eqref{LP-Z}.  We first exploit the compactness of $\Mbar$ to work on smaller regions.

\begin{proposition}[Local $L^4$ bound]\label{near-region}  
For every $z_0 \in M$ we have
$$
\int_0^1 \int_{B(z_0,\eta)} |u(z,t)|^4 \, dg(z) \, dt \leq C_{z_0,\eta} Z$$
if $\eta$ is sufficiently small (depending on $z_0$), where $B(z_0,\eta)$ is the ball of radius $\eta$
centered at $z_0$.
\end{proposition}

\begin{proposition}[Scattering $L^4$ bound]\label{far-region}  We have
$$
\int_0^1 \int_{\langle z \rangle > r_0} |u(z,t)|^4 \, dg(z) \, dt \leq 
C_{r_0} Z$$
if $r_0$ is sufficiently large.
\end{proposition}

From these two propositions and the compactness of $\overline{M}$ we easily
obtain \eqref{LP-Z} as desired.

It remains to prove the propositions.  Proposition \ref{near-region} is the
easier of the two and is proven in Section \ref{sec:near}, by applying the
positive commutator method on the product manifold $M \times M$ (with
Hamiltonian $H_{M \times M} = H_{z'} + H_{z''}$) to the product solution
$U(t,z',z'') := u(t,z') u(t,z'')$, and with the commutant $A_\near =
A_{\near,z_0,\eta}$ defined by
\begin{equation}
A_{\near} := \varphi(z',z'') i[H_{M \times M}, d_M(z', z'')],
\label{Anear}\end{equation}
where $\varphi(z',z'')$ is a smooth non-negative cutoff to the region $z', z'' \in B(z_0,\eta)$ which equals 1
on the region $z',z'' \in B(z_0,\eta/2)$.  Note that while the first-order
operator $A_\near$ is not quite self-adjoint, it does have real principal
symbol and its commutator will turn out to be positive definite (modulo lower order
terms).  The commutator $[H_{M \times M}, A_{\near}]$ turns out in fact to be quite tractable, 
as most of the error terms generated can be treated by the standard local smoothing estimates \eqref{eq:manifold-doi-comp}, \eqref{eq:manifold-doi-noderivs}.

Proposition \ref{far-region} is proven in a similar manner but involves
more technical difficulties (in particular, one needs much more precise
control on the derivatives of the metric function $d_M$ near infinity).
The proof applies the positive commutator method on $M \times M$ to the same
product function $U$ discussed before, but with the commutant $A_\asympt =
A_{\asympt,r_0}$ now given by
\begin{equation}
A_{\asympt} = \chi(z',z'') \psi(y', y'') i[H_{M \times M}, d_M(z', z'')]
\label{Aasympt}\end{equation} where $\chi(z',z'')$ is a smooth non-negative
cutoff to the product scattering region $\{ \langle z' \rangle, \langle z''
\rangle > r_0/2\}$ which equals 1 on the slightly smaller region
$\{ \langle z' \rangle, \langle z'' \rangle \geq r_0\}$, and $\psi$ is a
smooth non-negative cutoff\footnote{It may seem more intuitive to localize
$y'$ and $y''$ separately rather than to localize $d_{\partial M}(y',y'')$,
but as it turns out we cannot afford to let this cutoff be non-constant on
the diagonal $z' = z''$ and so we must proceed in this fashion.} to the
region $\{ d_{\partial M}(y',y'') \leq 2\eta\}$ which equals one on $\{
d_{\partial M}(y',y'') \leq \eta\}$, where $\eta > 0$ is a small number
(less than a quarter the injectivity radius of $\partial M$) to be chosen
later.  Note that while $\chi\psi$ localizes $y'$ and $y''$ to be close to
each other, it does not localized $r'$ and $r''$ to be similarly close
(although we have made them both large); it turns out that we cannot afford
to localize further as this will generate error terms that are too
difficult to estimate.

In the asymptotic regime, when controlling $[H, A_{\asympt}]$, terms arising from
derivatives hitting cutoff functions are harder to control, mainly because there are many terms which decay
only like $1/\langle z' \rangle$ or $1/\langle z''\rangle$, and the estimates \eqref{eq:manifold-doi-comp}, \eqref{eq:manifold-doi-noderivs} instead require decay of $1/\langle z'\rangle^{1+\eps}$, $1/\langle z''\rangle^{1+\eps}$.  One can use \eqref{eq:manifold-doi-ang} to control those terms which have two angular
derivatives.  Unfortunately there are also mixed terms involving
exactly \emph{one} angular derivative; these terms come from one derivative
hitting the angular cutoff $\psi$ and one hitting the
distance function $d_M$. A typical such term has the form
\begin{equation*}
\int_0^1 \int_{r', r'' \geq r_0/2}
\frac{a^{ij}(z', z'') \nabb_i u(t,z') \overline{\nabla_j u(t,z') }}{r'}
|u(t,z'')|^2 \, dg' \, dg'' \, dt,
\end{equation*}
where $a^{ij}(z', z'')$ are symbols of order 0 in $z'$, uniformly 
in $z''$ (see Proposition~\ref{cor:dist} below), and
$\nabb_i$ represents a component of the angular derivative $\frac{1}{r}
\nabla_y$. If both derivatives were angular then such a term could
be controlled by the single particle Morawetz estimate \eqref{eq:manifold-doi-ang},
however this estimate is insufficient when one of the derivatives has a
radial component.  Fortunately these mixed terms can still be handled by
the following new endpoint smoothing estimate:

\begin{lemma}\label{lemma:superDoi} Let $W$ be an arbitrary index set (possibly infinite).  For each $w \in W$, let
$a^{jk}_w(z)$ be a tensor-valued symbol ($j, k$ range over the co-ordinate indices of $\partial M$ and $M$ respectively) supported on the scattering region $\{ \langle z \rangle > r_0\}$ for sufficiently large $r_0$
which are uniformly symbols of order 0
in the sense that
$$ |\nabla_z^m a^{jk}_w(z)| \leq C_{m} \langle z \rangle^{-m}$$
for all\footnote{In fact, only finitely many $m$ are necessary.} $m \geq 0$, where the constants $C_m$ do not depend on
$w$.  Then we have
\begin{equation}\label{super-1}
\int_0^1 \sup_{w \in W} \bigg| \int_{M} \frac{a^{jk}_w(z) \nabb_j u(t,z)
\overline{\nabla_k u(t,z) }}{r} \, dg(z) \bigg| dt \leq C \| u(0)
\|_{H^{1/2}(M)}^2.
\end{equation} 
We also have the variant
\begin{equation}\label{super-2}
\int_0^1 \sup_{w \in W} \bigg| \int_{M} \frac{a^{j}_w(z) \nabb_j u(t,z)
\overline{u(t,z) }}{r} \, dg(z) \bigg| dt \leq C \| u(0)
\|_{L^2(M)}^2
\end{equation} 
where $a^j_w$ obeys similar bounds but without the additional index $k$.
\end{lemma}

This lemma appears to be new even when $W$ is a singleton set, and may be
of independent interest.  Morally speaking, it corresponds semi-classically
to the fact that $\int_{\langle z \rangle > r_0} |v_\ang(s)| / \langle z(s)
\rangle\ ds$ converges for all geodesics $z(s)$ of unit speed (cf.\ the
discussion after Lemma \ref{superDoi-intro}); however the result is still
rather subtle, since the inequalities fail if the absolute value signs are
brought inside the $z$-integral; see Proposition~\ref{counterex}.  The
proof is somewhat technical, requiring use of the scattering
pseudo-differential calculus, and is deferred to an appendix (Section
\ref{sec:doi}).  We remark in the case where $M$ is asymptotically flat
rather than asymptotically conic, one can avoid the need for this Lemma by
choosing a more sophisticated commutant (see the remarks in Section
\ref{remarks-sec}).

Note that in our applications of Lemma~\ref{lemma:superDoi}, the symbols
$a^{jk}_w(z)$ and $a^j_w(z)$ will often be derivatives of cutoff functions
and of the distance function on $M\times M,$ with $w$ corresponding to one
variable and $z$ to the other.

It remains to prove Proposition \ref{near-region} and Proposition \ref{far-region}.  This will be done in the next
few sections.  We first however record a number of quantities which are bounded by $Z$.

\begin{lemma}\label{usual-suspects} We have the estimates
\begin{align}
\label{us-1}
\int_0^1 \int_{M \times M} ( |\nabla u(t,z')|^2 |u(t,z'')|^2 +  |u(t,z')|^2 |\nabla u(t,z'')|^2 )& \\
\times (\frac{1}{\langle z' \rangle^2} + \frac{1}{\langle z'' \rangle^2})\, dg(z') \, dg(z'') dt &= O(Z)\nonumber\\
\label{us-1a}
\int_0^1 \int_{\langle z' \rangle, \langle z'' \rangle \geq r_0} |\nabb u(t,z')|^2 |u(t,z'')|^2 \frac{1}{r'}\, dg(z') \, dg(z'') dt &= O(Z)\\
\label{us-2}
\int_0^1 \int_{M \times M} |u(t,z')|^2 |u(t,z'')|^2 (1 + \frac{1}{d_M(z',z'')^2})\, dg(z') \, dg(z'') dt &= O(Z).
\\
\label{us-3}
\int_0^1 \int_{M \times M} \frac{a^{jk}(z',z'') \nabb_j u(t,z') \overline{\nabla_k u(t,z')}}{r'}
|u(t,z'')|^2\ dg(z') dg(z'')\ dt &= O(Z)\\
\label{us-4}
\int_0^1 \int_{M \times M} \frac{a^{jk}(z',z'') \nabb_j u(t,z') \overline{u(t,z')}}{r'} 
u(t,z'') \overline{\nabla_k u(t,z'')} dg(z') dg(z'')\ dt &= O(Z).
\end{align}
whenever $a^{jk}(z',z'')$ are symbols of order 0 in $z'$ uniformly in $z''$ that are
supported on the scattering region $\{ \langle z' \rangle\geq r_0 \}$ 
for $r_0$ sufficiently large.
\end{lemma}

\begin{proof}
We begin with \eqref{us-1}.  Consider for instance the term
$$
\int_0^1 \int_{M \times M} |\nabla u(z')|^2 |u(z'')|^2 \frac{1}{\langle z' \rangle^2}\, dg(z') \, dg(z'') dt.$$
Applying \eqref{eq:manifold-doi-comp} in the $z'$ variable and \eqref{spectrum} in the $z''$ variable
we can bound this expression by $O( \| u(0) \|_{H^{1/2}(M)}^2 \| u(0) \|_{L^2(M)}^2 )$,
which by \eqref{spectrum} again is $O( 2^k \| u(0) \|_{L^2(M)}^4 )$ $F= O(Z)$.  Now consider the term
$$
\int_0^1 \int_{M \times M} |u(z')|^2 |\nabla u(z'')|^2 \frac{1}{\langle z' \rangle^2}\, dg(z') \, dg(z'') dt.$$
Now we use \eqref{eq:manifold-doi-noderivs} in the $z'$ variable and \eqref{spectrum-2} in the $z''$ variable
to bound this by $O( \| u(0) \|_{H^{-1/2}(M)}^2 2^{2k} \| u(0) \|_{L^2(M)}^2 )$, which by \eqref{spectrum} is again
$O( 2^k \| u(0) \|_{L^2}^4 ) = O(Z)$.

The estimate \eqref{us-1a} is proven similarly to \eqref{us-1}, but with \eqref{eq:manifold-doi-ang} used
instead of \eqref{eq:manifold-doi-comp}.

Now consider \eqref{us-2}.  The integral on the region $d_M(z',z'') < \eps$ is clearly bounded by
$O( \eps^{-2} \| f \|_{L^2(M)}^4 ) = O(Z)$ thanks to \eqref{spectrum}.  When $d_M(z',z'') < \eps$ we observe that the kernel $(d_M(z',z''))^{-2}$ has marginal integrals
$$ \sup_{z''} \int_{d_M(z',z'') < \eps} (d_M(z',z''))^{-2} \ dg(z'), \quad
\sup_{z'} \int_{d_M(z',z'') < \eps} (d_M(z',z''))^{-2} \ dg(z'') \leq C \eps$$
so by Schur's test we have
\begin{align*}
 \int\int_{d_M(z',z'') < \eps} &|u(t,z')|^2 |u(t,z'')|^2 (d_M(z',z''))^{-2} \ dg(z') \,
dg(z'') \\
&\leq C \eps \| |u(t)|^2 \|_{L^2(M)} \| |u(t)|^2 \|_{L^2(M)}.
\end{align*}
Integrating this in $t$ we see that this contribution to \eqref{us-2} is
$O(\eps \int_0^1 \int_M |u(t,z)|^4\ dg dt)$ $= O(Z)$ as desired.

Now we prove \eqref{us-3}.  From \eqref{spectrum} we have
$$ \int |u(t,z'')|^2\ dg(z'') = O( \| u(0) \|_{L^2(M)}^2 )$$
so we can bound the left-hand side of \eqref{us-3} as
$$ O\big( \| u(0) \|_{L^2(M)}^2 \int_0^1 \sup_{z'' \in M} \big\lvert \int_M
\frac{a^{ij}(z',z'') \nabb_i u(t,z') \overline{\nabla_j u(t,z')}}{(r')^2}\
dg(z')\big\rvert dt \big )$$ which by Lemma \ref{lemma:superDoi} (letting the index
set $W$ be the manifold $M$) and \eqref{spectrum} is $O( \| u(0)
\|_{L^2(M)}^2 \| u(0) \|_{H^{1/2}(M)}^2 ) = O( 2^k \| u(0) \|_{L^2(M)}^4 )
= O(Z)$ as desired.

Finally we prove \eqref{us-4}.  From \eqref{spectrum}, \eqref{spectrum-2} and Cauchy-Schwarz we have
$$ \int |u(t,z'')| |\nabla u(t,z'')| \ dg(z'') = O( \| u(0) \|_{L^2(M)} \| u(0) \|_{H^1(M)} ) = O(2^k \| u(0) \|_{L^2(M)}^2 )$$
so by using Lemma \ref{lemma:superDoi} and \eqref{spectrum} as with \eqref{us-3} we can bound the left-hand side
of \eqref{us-4} as
$$ O( 2^k \| u(0) \|_{L^2(M)}^2 \| u(0) \|_{L^2(M)}^2 ) = O(Z)$$
as desired.
\end{proof}

\section{$L^4$ estimate in the near region}\label{sec:near}

In this section we prove the near region estimate Proposition \ref{near-region}, which is significantly
easier to prove (and indeed also follows from earlier work in \cite{st}, \cite{burq}) 
but already illustrates the basic method.  We fix $\eta$ and $z_0$, and allow our
constants to depend implicitly on $z_0$, $\eta$.

Let $A_\near$ be the commutant defined in \eqref{Anear}, and let $U = u
\otimes u$ as above.  From \eqref{adb}, Cauchy-Schwarz, \eqref{spectrum}, \eqref{spectrum-2} and the trivial observation
that $d_M$ is Lipschitz, we see that
$$\begin{gathered}
 |\langle A_\near U(t), U(t) \rangle_{M\times M}| \leq C \| U(t) \|_{H^1(M \times M)} \| U(t) \|_{L^2(M \times M)}\\
\leq C 2^k \| u(0) \|_{L^2(M)}^4 = O(Z)
\end{gathered}
$$
for $t=0,1$.  Thus by \eqref{heisenberg} we have
\begin{equation}
\label{near1}
\int_0^1\langle  i[H_{M \times M}, A_{\near}] U(t),  U(t) \rangle_{M \times M}\ dt \leq O(Z).
\end{equation}
The idea is now to decompose the left-hand side of \eqref{near1} as a positive main term that will yield Proposition
\ref{near-region}, plus some error terms which are either positive or $O(Z)$.

We first expand the commutator using \eqref{Anear} as
\begin{equation}\begin{gathered}
i[H_{M \times M}, \Anear] = -[H_{M \times M}, \varphi] [H_{M
\times M}, d_M(z', z'')] \\ - \varphi [H_{M \times M}, [H_{M
\times M}, d_M(z', z'')]].
\end{gathered}\label{near2}\end{equation}
Consider the contribution of the first term $[H_{M \times M}, \varphi] [H_{M
\times M}, d_M(z', z'')]$ in \eqref{near2} to \eqref{near1}.  This is an error term and will
be bounded in absolute value.  By integration by parts,
we can bound this contribution by
\begin{equation}\label{near2-first}
O(\int_0^1 |\langle  [H_{M \times M}, d_M] U(t),  [H_{M \times M}, \varphi] U(t) \rangle_{M \times M}|\ dt).
\end{equation}
To control this expression we apply \eqref{a-bracket}
to expand
\begin{equation}\label{hmdm}
 [H_{M \times M}, d_M(z', z'')] =  \langle \nabla_{z'} d_M, \nabla_{z'} \rangle_{g}
+  \langle (\nabla_{z''} d_M), \nabla_{z''} \rangle_{g}
+  (H_{z'} d_M) +  (H_{z''} d_M).
\end{equation}
To estimate these quantities, we first recall the form of the Laplacian in local 
polar co-ordinates:

\begin{lemma}\label{normal} Let $x^i$ be normal co-ordinates on a Riemannian manifold $X$, centred at $p \in X$, and let $(s, \phi)$ be polar co-ordinates with respect to $x^i$. That is, $s^2 = \sum_i (x^i)^2$ and $\phi$ are homogeneous of degree zero with respect to $x^i$. Then the Laplacian takes the form
$$ \Delta_X = - \frac{\partial^2 }{\partial s^2} - (1 + c) \frac{n-1}{s}
\frac{\partial}{\partial s} + s^{-2} \Delta_{S^{n-1}_s}.
$$ Here $c$ is a smooth function which is $O(s^2)$ and
$\Delta_{S^{n-1}_{s_0}}$ is the Laplacian on $S^{n-1}$ determined by the
co-ordinates $\phi$ and the metric $s_0^{-2} g$ restricted to the submanifold
$\{ s = s_0 \}$.
\end{lemma}

\begin{proof}The Laplacian in Riemannian polar normal co-ordinates is given by
$$
\frac1{\sqrt{g}} \,  \partial_s \sqrt{g} \, \partial_s + \frac1{\sqrt{g}} \partial_{\phi^i} \sqrt{g} \, g^{ij} \partial_{\phi^j}, \quad g = ds^2 + k_{ij}(s, \phi) d\phi^i d\phi^j, \quad k_{ij} = O(s^2). 
$$
Since in normal co-ordinates the metric is Euclidean to second order at the origin, $\sqrt{g} = s^{n-1}(1 + O(s^2))$. The result follows readily from these facts. 
\end{proof}

From this and \eqref{hmdm} we see that in the region $z',z'' \in B(z_0,\eta)$ we can write
$$
i[H_{M \times M}, d_M(z', z'')] = O(1) \nabla_{z',z''} + O(1) / d_M(z',z'')
$$
where $O(1)$ denotes various bounded functions of $z', z''$.  Applying this estimate, as well
as using \eqref{a-bracket} to expand $i[H_{M \times M}, \varphi]$, we can thus bound
\eqref{near2-first} by
\begin{align*}
 O(\int_0^1& \int_{B(z_0,2\eta) \times B(z_0,2\eta)}
(|U(t,z',z'')| + |\nabla_{z',z''} U(t,z',z'')|)\\
&
(|U(t,z',z'')|/d_M(z',z'') + |\nabla_{z',z''} U(t,z',z'')|)\ dg(z') dg(z'') dt.)
\end{align*}
Expanding out $U$ and using Cauchy-Schwarz, we can bound this by a linear combination of \eqref{us-1}, \eqref{us-2}.
Thus these terms are $O(Z)$.   

Comparing the above estimates with \eqref{near1} and \eqref{near2}, we now have
\begin{equation}\label{near-remaining}
\int_0^1 \langle -\varphi [H_{M \times M}, [H_{M \times M}, d_M]] U(t),  U(t) \rangle_{M \times M}\ dt \leq O(Z).
\end{equation}
Next, we apply \eqref{double-bracket} to obtain
\begin{equation}
-[H_{M \times M}, [H_{M \times M}, d_M]]
=
-\nabla_\beta \Hess(d_M)^{\alpha \beta} \nabla_\alpha
- H_{M \times M}^2 d_M
\label{near3}\end{equation}
where the gradients $\nabla$ are on the product manifold $M \times M$ (so the indices
$\alpha, \beta$ range over six values).

Now consider the contribution of the lower order term $H_{M \times M}^2 d_M(z', z'')$
to \eqref{near-remaining}.  We need

\begin{lemma}\label{compact-bidelta}  For $z', z'' \in B(z_0,2\eta)$, we have
$$ -H_{M \times M}^2 d_M(z',z'') = 8\pi \delta_{z'}(z'') + O(s^{-1}).
$$ 
\end{lemma}

\begin{proof} We again 
use polar normal
co-ordinates. Choose a co-ordinate system so that for each fixed $z''$ the
$z'$ co-ordinates are normal centred at $z''$. Then changing to polar
co-ordinates $s, \phi$ in the $z'$ variable we see that
\begin{equation*}\begin{gathered}
\Delta_{z'} d_M(z', z'') = \Big( - \frac{\partial^2 }{\partial s^2} - (1 + c)
\frac{n-1}{s} \frac{\partial}{\partial s} + s^{-2} \Delta_{S^{n-1}_r} \Big)
s \\ = (1 + c) \frac{n-1}{s},
\end{gathered}\end{equation*}
where $c(z',z'')$ vanishes to second order at $z' = z''$. This condition is
symmetric under interchanging $z'$ and $z''$ so the same is true of
$\Delta_{M \times M} d_M(z', z'')$; specializing to $n=3$, we thus have
$$ \Delta_{M \times M} d_M(z', z'') = (4 + \tilde c) d_M(z',z'')^{-1}$$
where $\tilde c$ also vanishes to second order at $z'=z''$.
Now we apply $\Delta_{M \times M}$ again using Lemma \ref{normal} to obtain
$$ \Delta_{M \times M}
\Delta_{M \times M} d_M(z', z'') = - 32 \pi \delta_{z'}(z'') + O(d_M(z',z'')^{-1})$$
and the claim follows (note that this is consistent with \eqref{flat-delta-delta}).
\end{proof}

Thus the contribution of this term to \eqref{near-remaining} is of the form
\begin{align*}
& 8\pi \int_0^1 \int_M \varphi(z,z) |u(t,z)|^4\ dg(z) 
\\&+ O( \int_0^1 \int_{B(z_0,2\eta) \times B(z_0,2\eta)} |u(z'')|^2 |u(z')|^2 / d_M(z,z')\ dg(z') dg(z'') ).
\end{align*}
The first term is the main term (as in \eqref{flat-delta-delta} and \eqref{interaction}).
The error term is $O(Z)$ by \eqref{us-2}.

The remaining term to consider in \eqref{near-remaining} is the contribution of $-\nabla \Hess(d_M) \nabla$,
which we can write using \eqref{inner-def} and integration by parts as
\begin{equation}\label{leib}
\begin{split}
&\Re \int_0^1 \int_{M \times M} \varphi \Hess(d_M)^{\alpha \beta} \nabla_\alpha U(t) \overline{\nabla_\beta U(t)}
\ dg(z') dg(z'') dt \\
&+
\Re \int_0^1 \int_{M \times M} (\nabla_\beta \varphi) \Hess(d_M)^{\alpha \beta} \nabla_\alpha U(t) \overline{U(t)}
\ dg(z') dg(z'') dt.
\end{split}
\end{equation}
The first term in \eqref{leib} is the most interesting one. We use the characterization of
$\Hess(d_M)(v,v)$, where $v=(v',v'') \in T_{(z',z'')}M^2$, as the second derivative of
$d_M(z'(t), z''(t))$ where $(z'(t), z''(t))$ moves along a geodesic in $M^2$
with initial condition $((z', z''); v)$. We use the second variation formula for 
geodesics (see \cite{jost}, Theorem 4.1.1) to obtain
\begin{equation}\begin{gathered}
\frac{d^2}{dt^2} d_M(z'(t), z''(t)) \geq 
- \frac1{\dist(z',z'')} \int_0^1 \angs{R(\pa_s,\pa_t^\perp)\pa_t^\perp, \pa_s} \, ds  \\ \geq - C d_M(z'(t),z''(t)) \big(|v'|^2 + |v''|^2 \big),
\end{gathered}\label{svf}\end{equation}
where $\partial_s$ is the velocity of the geodesic which goes from $z'(0)$
to $z''(0)$ in unit time, hence $|\partial_s| = d_M(z',z'')$ and
$\partial_t$ is the Jacobi field corresponding to the variation of the
geodesic.  In particular this bound does not blow up when $d_M(z',z'')$ is
small. From this it follows that
\begin{align*}
\Re \int_0^1 \int_{M \times M} &\varphi \Hess(d_M)^{\alpha \beta} \nabla_\alpha U(t) \overline{\nabla_\beta U(t)}
\ dg(z') dg(z'') dt \\
&
\geq -C \int_0^1 \int_{B(z_0,2\eta) \times B(z_0,2\eta)} |\nabla U(t)|^2\ dg(z') dg(z'') dt.
\end{align*}
The right-hand side is then $O(Z)$ by \eqref{us-1}.

Finally, consider the second term in \eqref{leib}.  Using the crude estimate $|\Hess(d_M)| \leq C / d_M$
we can bound this contribution by
$$ O( \int_0^1 \int_{B(z_0,2\eta) \times B(z_0,2\eta)} |U(t,z',z'')| |\nabla_{z',z''} U(t,z',z'')| / d_M(z',z'')\ dg(z') dg(z'') dt).$$
By Cauchy-Schwarz we can bound this a linear combination of \eqref{us-1} and \eqref{us-2}, and so these
terms are also $O(Z)$.

Inserting all of the above estimates into \eqref{near-remaining}, we obtain Proposition \ref{near-region}.


\section{The geometry of the scattering region}

To prove Theorem \ref{main}, it remains to verify the scattering region $L^4$ estimate in 
Proposition \ref{far-region}.  To do this we shall need to understand the geometry of the scattering
region.  More precisely, we will
need to control the metric function $d_M(z',z'') = d_M((x',y'),(x'',y''))$ and its derivatives on the support
of $\chi\psi$, where $\chi$, $\psi$ are as in \eqref{Aasympt}.
To avoid singularities forming in $d_M$, we will always assume in this section
that $r_0$ is sufficiently large and $\eta$ sufficiently small.

We begin with a basic symbol estimate on the metric $d_M$.

\begin{proposition}\label{symbol-est}  Let $(z', z'')$ be in the support of $\chi\psi$.  If $d_M(z',z'') \geq \frac{r'+r''}{100}$, then we have the product symbol estimates
$$ |\nabla_{z'}^{m'} \nabla_{z''}^{m''} \nabla_{z',z''} d_M(z',z'')| \leq C_{m',m''} 
(r')^{-m'} (r'')^{-m''}$$
for all $m',m'' \geq 0$.   If instead $d_M(z',z'') \leq \frac{r'+r''}{100}$, then we have
$$ |\nabla_{z',z''}^m d_M(z',z'')| \leq C_m d_M(z',z'')^{1-m}$$
for all $m \geq 0$.  Furthermore, for $j \leq 3$, the distribution
$\nabla_{z',z''}^j d_M(z',z'')$ has no mass on the diagonal $z'=z''$.
\end{proposition}
This symbol estimate, while quite plausible, is a little technical and its proof will be
deferred to Section \ref{appendix}.  From this, 
Lemma \ref{compact-bidelta}, and a rescaling argument we can obtain the following estimate for the 
bi-Laplacian of $d_M$:

\begin{corollary}\label{lemma:lapsquareddist}
The function $H_{M \times M}^2 d_M(z',z'')$ satisfies the distributional inequality
$$
-H_{M \times M}^2 d_M(z',z'') \leq 8 \pi \delta_{z'}(z'') - C d_M(z',z'')^{-1} (r' + r'')^{-2}
- C (r')^{-3} - C (r'')^{-3}
$$
uniformly on the support of $\chi\psi$.
\end{corollary}

\begin{proof} Observe that if $|r''/r'-1| > c$ for some $c>0$ $d_M(z',z'')$ is comparable to $r'+r''$ and
the estimate follows from Proposition \ref{symbol-est}, so we may assume that $|r''/r'-1| < c$ for some small $c$.  
In this region Proposition \ref{symbol-est} is insufficient (it gives a bound of $O(d_M^{-3})$, which is not
even locally integrable near the diagonal).  Instead, we
we set $\tau := 1/r'$ and consider the rescaled manifold\footnote{Strictly speaking, this rescaling is
only well-defined on the scattering region of the manifold, but that is all we shall need here.}
 $M_\tau$ with metric
$$
g_\tau = dr^2 + r^2 h_{ij}(\frac{\tau}{r},y) dy^i dy^j.
$$ Observe that $g_\tau$ various smoothly in $\tau$ as $\tau$ varies in the
compact set $[0, 2/r_0]$.  Also, since $d_{\partial M}(y',y'') \leq 2\eta$
and $|r''/r' - 1| < c$, we see that
$$d_{M_\tau}(\tau z',\tau z'') \leq C(c + \eta),$$
with the action $\tau z' := (\tau r', y')$ and likewise in $z''.$
Thus if $c$ and $\eta$ are sufficiently small, we can apply Lemma \ref{compact-bidelta} and conclude that
$$ -H_{M_\tau \times M_\tau}^2 d_{M_\tau}(z',z'') \leq 8 \pi \delta_{\tau z'}(\tau z'') - C d_{M_\tau}(\tau z',\tau z'')^{-1}$$
uniformly in $\tau$.  The claim follows by undoing the scaling.
\end{proof}
  
Finally, we need a geodesic convexity estimate on the metric $d_M(z'(s), z''(s))$, generalizing \eqref{psd}.

\begin{proposition}\label{full-ang-est} Let $z'(s), z''(s)$ be two geodesics.
For $r_0$ sufficiently large and $\eta$ sufficiently small, the estimate 
\begin{equation}\begin{gathered}
\frac{d^2}{ds^2} d_{M}(z'(s), z''(s)) \geq - C\Big( \frac{|v'_{\ang}|^2}{r'} + \frac{|v''_{\ang}|^2}{r''} \Big) \\ 
- C \Big( \frac{|v'|^2}{(r')^2} + \frac{|v''|^2}{(r'')^2} \Big) 
\end{gathered}\label{ang-est2}\end{equation}
holds whenever $z'(s), z''(s)$ are distinct points in the support of $\chi\psi$, 
where $|v'| := |\frac{d}{ds} z'(s)|_{g(z'(s))}$ is the speed of $z'$ and $|v'_\ang| := r'(s) |\frac{d}{ds} y'(s)|_{h(z'(s))}$ is the angular component of this speed, and similarly for $v''$.
\end{proposition}

This proposition is somewhat delicate and will be also deferred to Section
\ref{appendix}.  The idea is to first prove this lemma for perfectly conic
metrics (in which the second group of terms in \eqref{ang-est2} do not
appear); the perfectly conic metric contains a large number of totally
geodesic planes on which the metric behaves like the Euclidean metric
except in normal directions (generating the first group of terms in
\eqref{ang-est2}).  We then handle the approximately conic case by
obtaining symbol estimates for the error between the approximately conic
and perfectly conic metrics giving rise to the second group of terms in
\eqref{ang-est2}.

\section{$L^4$ estimate in the scattering region}\label{sec:asympt}

We can now begin the proof of Proposition \ref{far-region}.  Fix $r_0, \eta$; we assume $r_0$ to be large enough
and $\eta$ small enough that the geometrical lemmas of the previous section apply, and allow our constants $C$ to depend on $r_0$, $\eta$.  
Let $\Aasympt$ be
as in \eqref{Aasympt}.  By arguing as in Section \ref{sec:near}
we have
\begin{equation}
\label{asympt1}
\int_0^1\langle  i[H_{M \times M}, A_{\asympt}] U(t),  U(t) \rangle_{M \times M}\ dt \leq O(Z).
\end{equation}
Again, we will expand the left-hand side of \eqref{asympt1} to extract a positive main term that
will give the estimate in Proposition \ref{far-region}, plus a number of error terms which are either positive or
$O(Z)$.

From \eqref{Aasympt} we may decompose
\begin{equation}\begin{gathered}
i[H_{M \times M}, \Aasympt] = -[H_{M \times M}, \chi]
\psi [H_{M \times M}, d_M] \\ - \chi [H_{M
    \times M}, \psi] [H_{M \times M}, d_M]\\ - \chi \psi [H_{M \times M}, [H_{M \times M}, d_M]].
\end{gathered}
\label{asympt2}\end{equation}
The first two terms are error terms and will be estimated in absolute value; it is the third one which will generate the interesting positive terms.

We first analyze the first term $[H_{M \times M}, \chi] \psi [H_{M \times M}, d_M]$ in \eqref{asympt2}.  Integrating by parts we see that the contribution of this term to \eqref{asympt1} is
\begin{equation}\label{asympt2-first}
 O( \int_0^1 |\langle \psi i[H_{M \times M}, \chi] U(t), 
i[H_{M \times M}, d_M] U(t) \rangle_{M \times M}|\ dt).
\end{equation}
This integral is supported in the region
where either $\langle z'\rangle \leq r_0$ or $\langle z'' \rangle \leq r_0$, since $\chi$ is constant elsewhere; without loss of generality it suffices to consider
the terms when $\langle z'\rangle \leq r_0$.  To control this expression, we first observe from \eqref{a-bracket}
and Lemma \ref{symbol-est} that
\begin{equation}\label{lap-est} i[H_{M \times M}, d_M] = O(1) \nabla_{z',z''} + O(1) / d_M(z',z'')
+ O(1)
\end{equation}
on the support of $\chi\psi$, where we use $O(1)$ to denote various bounded (tensor-valued) functions of $z', z''$;
for future reference we also remark that these $O(1)$ errors obey symbol estimates when $d_{\partial M}(y',y'') \geq \eta/2$, since in this case $d(z',z'') \geq c_\eta (r' + r'')$.

From \eqref{lap-est}, \eqref{a-bracket} and Cauchy-Schwarz, we thus see that \eqref{asympt2-first} can be
bounded in magnitude by
\begin{align*}
 O( \int_0^1 \int_{\langle z' \rangle \leq r_0}&
|\nabla u(t,z')|^2 |u(t,z'')|^2 +
|u(t,z')|^2 |\nabla u(t,z'')|^2 \\
&+ |u(t,z')|^2 |u(t,z'')|^2 (1 +  d_M(z',z'')^{-2})\ dg(z') dg(z'') ).
\end{align*}
The first two terms are $O(Z)$ by \eqref{us-1}.  For the third term, observe that the contribution of
the region $r'' < 4r_0$ is $O(Z)$ by \eqref{us-2}, while for the region $r'' > 4r_0$ we simply
estimate $(1 + d_M(z',z'')^{-2})$ crudely by $O(1)$, and this term is then bounded using \eqref{spectrum}
by $O( \|u(0) \|_{L^2(M)}^4 ) = O(Z)$.  Thus the total contribution of the first term of \eqref{asympt2} to \eqref{asympt1} is $O(Z)$.

Now consider the contribution of the second summand
of \eqref{asympt2} to \eqref{asympt1}; these are also error terms, but of a more delicate nature.  Again
we integrate by parts to estimate this contribution by
\begin{equation}\label{asympt2-second}
 O( \int_0^1 |\langle i[H_{M \times M},\psi] \chi U(t), 
i[H_{M \times M}, d_M] U(t) \rangle_{M \times M}|\ dt).
\end{equation}
The key point here is that since the angular cutoff
$\psi$ does
not depend on the radial variables $r', r''$, the first-order operator
$[H_{M \times M}, \psi]$ only contains angular derivatives $\nabla_{y',y''}$ and
no radial derivatives $\nabla_{r',r''}$, thanks to \eqref{a-bracket}.  Indeed, on the support of
$\chi$ we have
\begin{align*}
 [H_{M \times M}, \psi] = &O(1) \frac{1}{(r')^2} \nabla_{y'}
+ O(1) \frac{1}{(r'')^2} \nabla_{y''} \\
&+ \frac{1}{(r')^2} O(1) + \frac{1}{(r'')^2} O(1)
\end{align*}
where again we use $O(1)$ to denote various bounded functions of $z', z''$ which also obey symbol
estimates.  From this, \eqref{lap-est}, 
 and symmetry we see that these contributions are bounded by a combination of
expressions of the form
\begin{align*}
O\bigg( \big|\int_0^1 \int_{M \times M} &\frac{a^{jk}(z',z'') \nabb_j u(t,z') \overline{\nabla_k u(t,z')}}{r'}
|u(t,z'')|^2\ dg(z') dg(z'')\ dt \big| \bigg),\\
O\bigg( \big|\int_0^1 \int_{M \times M} &\frac{a^{jk}(z',z'') \nabb_j u(t,z') \overline{u(t,z')}}{r'} 
u(t,z'') \overline{\nabla_k u(t,z'')} dg(z') dg(z'')\ dt \big| \bigg),\\
O\bigg( \big|\int_0^1 \int_{r', r'' \geq r_0} &\frac{|\nabla u(t,z')| |u(t,z')| |u(t,z'')|^2}{r' d_M(z',z'')}\ dg(z') dg(z'')\ dt \big| \bigg),\\
O\bigg( \big|\int_0^1 \int_{r', r'' \geq r_0} &\frac{|\nabla u(t,z')| |u(t,z')| |u(t,z'')|^2}{(r')^2}\ dg(z') dg(z'')\ dt \big| \bigg),\\
O\bigg( \big|\int_0^1 \int_{r', r'' \geq r_0} &\frac{|u(t,z')|^2 |u(t,z'')|^2}{(r')^2 d_M(z',z'')}\ dg(z') dg(z'')\ dt \big| \bigg)
\end{align*}
(noting that $\nabla_y = r \nabb$), where $a^{jk}$ denotes symbols
supported on the the support of $\chi \nabla \psi$ (note that we do not
encounter the singularity on the diagonal $z'=z''$ because $\nabla \psi$
vanishes here).  The first two terms are controlled by \eqref{us-3},
\eqref{us-4} respectively.  The third and fourth terms can be controlled
via Cauchy-Schwarz by a linear combination of \eqref{us-1} and
\eqref{us-2}.  Finally, the last term is easily controlled by \eqref{us-2}.
(Recall that $r'$ is large.)  Thus the total contribution of the second term of
\eqref{asympt2} to \eqref{asympt1} is $O(Z)$.

In light of \eqref{asympt1}, \eqref{asympt2}, and the preceding estimates, we see that
\begin{equation}\label{asympt-final}
\int_0^1\langle -\chi \psi [H_{M \times M}, [H_{M \times M}, d_M(z', z'')]] U(t),  U(t) 
\rangle_{M \times M}\ dt \leq O(Z).
\end{equation}
Of course, we may expand $-[H_{M \times M}, [H_{M \times M}, d_M(z', z'')]]$ using \eqref{double-bracket} as 
\begin{equation}\begin{gathered}
-[H_{M \times M}, [H_{M \times M}, d_M(z', z'')]]=-\nabla_\beta \Hess(d_M)^{\alpha \beta}
\nabla_\alpha - H_{M \times M}^2 d_M
\end{gathered}\label{asympt3}\end{equation}
where $\Hess(d_M)^{\alpha \beta} = \nabla^\alpha \nabla^\beta d$.

Let us first study the lower order term $-H_{M \times M}^2 d_M(z', z'')$.
Applying Corollary \ref{lemma:lapsquareddist}, the contribution of this
term to \eqref{asympt-final} is greater than or equal to
$$
\begin{gathered}
 \geq 8 \pi \int_0^1 \int_M \chi(z,z)^2 |u(t,z)|^4 \, dg(z) \, dt\\
+ O( \int_0^1 \int_{M \times M} \frac{|u(t,z')|^2 |u(t,z'')|^2}{d_M(z',z'') (r'+r'')^{2}} \ dg(z') dg(z'') dt).
\end{gathered}
$$ 
The error term here is of course $O(Z)$ by \eqref{us-2}.

It remains to treat the contribution of the $-\nabla \Hess(d_M) \nabla$ to \eqref{asympt-final}.  
Again we integrate by parts, creating a main term
\begin{equation}\label{main-term}
\Re \int_0^1 \int_{M \times M} \langle \chi(z',z'') \psi(y',y'') \Hess(d_M)^{\alpha \beta} 
\nabla_\alpha U(t,z',z'') \overline{\nabla_\beta U(t,z',z'')}\  dt
\end{equation}
and an error term which is bounded by 
$$ \int_0^1 \int\limits_{M \times M}   |\nabla_{z',z''} 
(\chi(z',z'') \psi(y',y'')| |\Hess(d_M)|
|\nabla_{z',z''} U(t,z',z'')| |U(t,z',z'')|\  dt.$$
Consider the error term first.  For this term we use the crude estimate
$$ |\Hess(d_M)(z',z'')| \leq C (1 + d_M(z',z'')^{-1})$$
from Lemma \ref{symbol-est}, as well as 
$$ |\nabla_{z',z''} (\chi(z',z'') \psi(y',y'')| \leq \frac{1}{r'} + \frac{1}{r''}.$$
By symmetry we can thus estimate this error term by
$$ C \int_0^1 \int_{M \times M} |\nabla u(t,z')| |u(t,z')| |u(t,z'')|^2 (\frac{1}{r'} + \frac{1}{r''}) (1 + \frac{1}{d_M(z',z'')})
\ dg(z') dg(z'') dt.$$
By Cauchy-Schwarz one can estimate this by a linear combination of \eqref{us-1} and \eqref{us-2}.

Now consider the main term \eqref{main-term}.  Recall that $\Hess(d_M)(v,v)$,
where $v \in T_{(z',z'')}M^2$, is the second derivative of $d_M(z'(s),
z''(s))$ in $s$ as $(z'(s), z''(s))$ moves along a geodesic in $M^2$ with initial
condition $((z', z''); v)$.
By Proposition \ref{full-ang-est}, we see that 
\begin{align*}
\eqref{main-term} \geq &
-C \int_0^1 \int_{r', r'' \geq r_0/2} 
\Big( \frac{|\nabb u(t,z')|^2 |u(t,z'')|^2}{r'} + \frac{|u(t,z')|^2 |\nabb u(t,z'')|^2}{r''} \Big)
 \\
&+
\Big( \frac{|\nabla u(t,z')|^2 |u(t,z'')|^2}{(r')^2} + \frac{|u(t,z')|^2 |\nabla u(t,z'')|^2}{
(r'')^2} \Big)
\ dg(z') dg(z'') dt.
\end{align*}
But the first term on the right-hand side is $O(Z)$ by \eqref{us-1a}, while the second term is also $O(Z)$ by
\eqref{us-1}.  Combining all of these estimates together we obtain Proposition \ref{far-region} (once we
verify Proposition \ref{cor:dist} and Proposition \ref{full-ang-est} in the next section).  This completes
the proof of Theorem \ref{main}.


\section{Concluding remarks}\label{remarks-sec}

We first want to comment on some particular geometric settings in which the
proof of the Morawetz estimate is much easier. The first is that of a simply-connected
smooth three-manifold $M$ which has globally non-positive, bounded
sectional curvature; by the Hadamard-Cartan theorem (see \cite{jost},
Corollary 4.8.1), such a manifold is diffeomorphic to $\R^3.$ Such a
manifold is automatically non-trapping, and the distance function
$d_M(z',z'')$ is smooth on all of $M \times M$ outside the diagonal.  We
also make the technical assumption that the distance function $d_M(z',z'')$
is uniformly $C^4$ in the off-diagonal region $\{ (z',z'') \in M \times M:
d_M(z',z'') > 1 \}$\footnote{It may be that this is automatically true, but
this is not clear to us.}.  In particular $H^2 d_M$ is bounded in this
region. More importantly the distance function is globally convex (see \cite{jost}, Theorem 4.1.1), so that
the tensor $\Hess(d_M)^{\alpha \beta}$ is positive semidefinite, by the
second variation formula. Thus with the commutant $A = i[H, d_M(z',z'')]$
(with no cutoffs whatsoever) we obtain a Morawetz inequality very easily
for this class of manifolds. 

The next case we consider is that of asymptotically flat three-dimensional
non-trapping manifolds, i.e.\ where $M = \R^3$ and the metric $g$ satisfies
\eqref{ac-metric} with $h$ equal to the standard metric on $S^2$.  We make
a distinction between the distance $d(z',z'')$ on $M$ induced by the metric
$g$ and the Euclidean distance $|z'-z''|$ obtained by identifying $M$ with
Euclidean space.  We shall only discuss the Morawetz estimate in the
scattering region $\langle z' \rangle, \langle z'' \rangle \geq r_0$.  Here
the key observation is that when $z'$ and $z''$ are fairly close together,
e.g.\ $|z'-z''| < \frac{1}{2} |z'|$, then the distance $d_M(z',z'')$ is
smooth and is close to $|z'-z''|$ in the sense that
$$ \bigl| d_M(z',z'') - |z'-z''|\bigr| \leq C \hbox{ whenever } 
\langle z' \rangle, \langle z'' \rangle \geq r_0
\hbox{ and } |z'-z''| < \frac{1}{2} |z'|,$$ at least if we
make the compact region $K_0$ sufficiently large.  This is an easy
consequence of the decay of the difference $g-\delta$ between the two
metrics.  In fact, the error obeys symbol estimates of one order better than
Proposition \ref{symbol-est} suggests; see Proposition \ref{cor:dist} in the appendix.
We now consider the commutant $A := i[H,a]$, where $a$ is the function
\begin{equation}\label{a-flat}
a := \chi(z',z'') (\varphi(z',z'') d_M(z',z'') + (1 - \varphi(z',z'')) |z'-z''|),
\end{equation}
$\varphi$ is a smooth cutoff which
equals 1 when $|z'-z''| \leq \frac{1}{4} |z'|$ and equals 0 when $|z'-z''|
\geq \frac{1}{2} |z'|$, and $\chi$ is as in Section \ref{overview}.  
Thus this cutoff is equal to the actual distance
$d(z',z'')$ when $z'$ and $z''$ are close, but reverts smoothly to the
Euclidean distance when $z'$ and $z''$ are far apart.  This significantly
reduces the error caused by differentiating the cutoff $\varphi$, indeed it
generates terms that are now controlled by \eqref{us-1} rather than the
terms \eqref{us-2}, \eqref{us-3}, \eqref{us-4} which required the
Morawetz estimate \eqref{ls-ang} and its refinement in 
Lemma \ref{lemma:superDoi}.  We omit the details.

Next, we give an example to show that \eqref{super-1} cannot be
strengthened by bringing the absolute value signs inside the $z$-integral.

\begin{proposition}\label{counterex}  There does not exist an estimate of the form
\begin{equation}\label{false}
\int_0^1 \int_{M} \frac{ \big| \nabb u(t,z) \big| \, \big| \nabla u(t,z)
\big| }{r} \, dg \, dt \leq C \| u(0) \|_{H^{1/2}(M)}^2 \quad \text{
(FALSE!)}
\end{equation}
for solutions to \eqref{nls} even in Euclidean space $(M,g) = (\R^3,
\delta)$.
\end{proposition}

\begin{proof}[Proof of falsity of \eqref{false}] We give a sketch only.
Pick a large integer $N = 2^k > 1$, and let $\phi$ be the explicit solution
to \eqref{nls} given by
$$ \phi(t,z) := (t+i)^{-3/2} e^{\frac{i|z|^2}{2(t+i)}};$$ note that at time
$t=0$ this is basically a standard Gaussian, while for times $t \gg 1$, the
solution $\phi(t)$ is concentrated on a ball of radius $O(t)$ centered at
the origin, and it (and its derivatives) have magnitude roughly $t^{-3/2}$
on this ball.  We then let $u$ be the function
$$ u(t,z) := k^{1/2} \phi(N^2 t, Nz) + \sum_{j=1}^k \epsilon_j \phi(N^2 t,
N(z - 2^j e_1))$$ where $e_1 = (1,0,0)$ is a standard unit vector and
$\epsilon_j = \pm 1$ are random signs chosen independently; the purpose of
the signs is to ensure (thanks to Khinchin's inequality) that we do not
expect any unusual cancellation between the summands.

The function $u(t,z)$ is clearly a solution to \eqref{nls}, and the
$H^{1/2}(\R^3)$ norm can be computed as
$$ \| u(0) \|_{H^{1/2}(\R^3)} \sim k^{1/2} N^{-1};$$ this is easiest shown
by first verifying the more general statement $\| u(0) \|_{H^s(\R^3)} \sim
k^{1/2} N^{s-3/2}$ for $s=0,1,2$ (applying a rescaling by $N$ if desired)
and then interpolating.  Thus the right-hand side of \eqref{false} is $O(k
N^{-2})$.

Let us choose any $N/3 \leq j \leq 2N/3$ and consider the size of $u(t,z)$
and its derivatives in the region of spacetime where $t \sim 2^j/N$ and $|z
- 2^j e_1| \ll 2^j$.  In this region $|\nabla k^{1/2} \phi(N^2 t, Nz)|$ has
size $\sim k^{1/2} (2^j N)^{-3/2}$, and so the expectation of $|\nabla
u(t,z)|$ is also at least $\gtrsim k^{1/2} (2^j N)^{-3/2}$.  Since $\phi$
is radial, the angular derivative $|\nabb k^{1/2} \phi(N^2 t, Nz)|$
vanishes, however $|\nabb \phi(N^2 t, N(z - 2^j e_1))|$ is fairly large,
comparable to $\sim (2^j N)^{-3/2}$ in a large fraction of the region under
consideration.  Thus this region of spacetime contributes at least $\gtrsim
k^{1/2} (2^j N)^{-3} (2^j N) 2^{3j} / 2^j \sim k^{1/2} N^{-2}$ to the
integral in \eqref{false}.  Summing over all $j$ we see that the left-hand
side of \eqref{false} is at least $\sim k^{3/2} N^{-2}$.  Comparing this
with the right-hand side we obtain the desired contradiction\footnote{Note
that this is only a `logarithmic' failure in the frequency variable.
Indeed if we concede an epsilon worth of derivatives then we can handle
these terms easily by a variant of \eqref{eq:manifold-doi-comp}
where the epsilon loss is transferred from the weight $r^{-1-\eps}$ to the
regularity $H^{1/2+\eps}(M)$, thus recovering the result of Burq
\cite{burq}, at least in the context of $L^4_{t,x}$ estimates.} by setting
$k \to \infty$.
\end{proof}

\begin{remark} One can construct a similar example to show that
the second estimate in Lemma \ref{lemma:superDoi} similarly fails if
the absolute values are placed inside the integral.
\end{remark}

\begin{remark}
Note that in our example $\nabb u(t,x)$ and (the largest term in) $\nabla
u(t,x)$ oscillate in different directions and so if one removes the
absolute values, replacing instead by a smooth symbol, then we do not
obtain a counterexample to Lemma \ref{lemma:superDoi}.
\end{remark}

Finally, we remark on the applicability of the Strichartz estimate in
Theorem \ref{main} to non-linear Schr\"odinger equations.  It turns out that
in such equations, the most useful Strichartz estimates are either those that
require no derivatives whatsoever on $u$ (i.e.\ they have $L^2(M)$ on the right-hand side),
or measure $u$ in a space of the form $L^q_t L^\infty_x$.  Our current estimate has neither
of these two properties.  However by commuting the Schr\"odinger flow with a power of $(1+H)$
and using Sobolev embedding one can get an $L^q_t L^\infty_x$ estimate, for instance we have
$$ \| u \|_{L^4_t L^\infty_x([0,1] \times M)} \leq C_\eps \| u(0) \|_{H^{1+\eps}(M)}$$
for any $\eps > 0$.  This result can for instance be used to demonstrate local well-posedness in $H^{1+\eps}(M)$ 
of the quintic non-linear Schr\"odinger equation $iu_t + \frac{1}{2} \Delta u = \pm |u|^4 u$ on non-trapping
asymptotically conic manifolds by a standard argument which we omit (see e.g.\ \cite{cwI}, \cite{bgt}).
Note that the above estimate was also derived by Burq \cite{burq} for asymptotically flat non-trapping manifolds.

It would be interesting to see if the $\eps$ can be removed; this could
then be used to demonstrate \emph{global} well-posedness of the above
quintic equation in the energy space $H^1(M)$ for small energy data, using
the conservation of energy in the standard manner (note the sign of the
non-linearity is irrelevant in the small energy setting).  While the $\eps$
can indeed be removed from the above Strichartz estimate in the Euclidean
space case (see \cite{tao:focusing}; the corollary concerning global
well-posedness can be achieved by many other means, see for instance
\cite{cwI}), we were unable to remove it here.  However, by a Besov space
interpolation argument (see \cite{tao:focusing}) it is possible to obtain
the variant estimate
$$ \| u \|_{L^q_t L^\infty_x([0,1] \times M)} \leq C_q \| u(0) \|_{H^{3/2 -
2/(q-1)}(M)}$$ for all $4 < q < \infty$.  This suffices to give local
well-posedness on non-linear Schr\"odinger equations $iu_t + \frac{1}{2}
\Delta u = F(u)$ in the space $H^{3/2 - 2/(q-1)}(M)$ for any $q > 4$, whenever
$F(u)$ is smooth and grows like $|u|^q$; we omit the details as they are
rather standard (see e.g.\ \cite{cwI}).  Note that in Euclidean space the
space $H^{3/2-2/(q-1)}$ would be critical with respect to the usual scaling of
the non-linear Schr\"odinger equation.  Thus the Strichartz estimate
Theorem \ref{main} has some application to non-linear Schr\"odinger
equations, though clearly it is less satisfactory than the full range
\eqref{cond} of Strichartz estimates would be; these would yield local
well-posedness in $H^s$ for any $q>1$ and $s \geq
\max(0, \frac{3}{2} - \frac{2}{q-1}).$


\section{Appendix I: Proof of geometrical lemmas}\label{appendix}

In this section we analyze the behaviour of geodesics in thin conic regions of $M$, in order
to prove the geometric estimates in Proposition \ref{symbol-est} and Proposition \ref{full-ang-est} required
in the above proof of Theorem \ref{main}.  The main strategy is to compare the asymptotically conic metric
$d_M$ to a perfectly conic metric $d_\conic$, which has an exact formula and can be handled explicitly.  

By compactness of $\partial M$ it suffices to work in a truncated sector 
\begin{equation}
\Omega(y_0, \eta, r_0) := \{ (r, y) \mid d_{\partial M}(y, y_0) < \eta; \langle z \rangle > r_0 \}
\label{Omega}\end{equation}
for some $y_0 \in \partial M$, where we assume $r_0$ sufficiently large and $\eta$ sufficiently small 
(possibly depending on $y_0$, although this is irrelevant since $y_0$ ranges in a compact set).

We first briefly review geodesic flow on a general Riemannian manifold $(M,g)$ parameterized by a general set of 
co-ordinates $z^j$.
We will view geodesic flow as a Hamiltonian flow on the cotangent bundle $T^*(M)$ of $M$, 
which we parameterize as $(z^j,\zeta_j)$, where $\zeta_j$ is the co-ordinate dual to the vector field $d z^j$; thus
$\zeta_j d z^j$ is the canonical one-form.  On this cotangent bundle we
define the energy function $\sigma(H): T^*(M) \to \R$
as
$$ \sigma(H) := \frac{1}{2} g^{jk}(z) \zeta_j \zeta_k = \frac{1}{2}
|\zeta|_{g(z)}^2;$$ note that this is also the symbol of the Schr\"odinger
operator $H$.  This energy function then induces trajectories $(z(s),
\zeta(s))$ in the cotangent bundle which are just the paths of geodesic
flow (using the metric $g$ to identify the cotangent and tangent bundles):
\begin{align*}
\frac{d}{ds} z^j &= g^{jk}(z) \zeta_k \\
\frac{d}{ds} \zeta^j &= - \frac{1}{2} \frac{\partial g^{kl}}{\partial z^j}(z) \zeta_k \zeta_l.
\end{align*}
These flows preserve $\sigma(H)$, and are also homogeneous
with respect to the scaling $(z(s), \zeta(s)) \mapsto (z(\lambda s),
\lambda \zeta(\lambda s))$, which corresponds to speeding up the velocity
of the geodesic flow by $\lambda$.  If we use this homogeneity to normalize
to the unit speed case $\sigma(H) = \frac{1}{2}$ (i.e.\ $|\zeta|_{g(z)} =
1$), then this flow is related to the distance function $d_M$ on $M$ by the
formula
$$ d_M( z(s), z(s') ) = |s-s'|$$
provided that $|s-s'|$ is sufficiently small.  Thus the short-time geodesic flow can be used to define the metric locally.

Now we specialize to studying geodesic flow in the truncated sector $\Omega(y_0, \eta, r_0)$ in the scattering
co-ordinates $(x, y^j)$, where the manifold $M$ is either asymptotically conic or perfectly conic.  
Of course this region is not geodesically complete, so the statements below are
conditional assuming that the geodesic does indeed stay inside this region.

The first step is to choose co-ordinates for the cotangent bundle $T^*
\Omega(y_0, \eta, r_0)$.  We shall use the co-ordinates
$(x,y^j,\nu,\mu_j)$, where $0 < x < \epsilon_0$, and $y \in \partial M$,
while the co-ordinates $\nu$ and $\mu_j$ are dual to the vector fields
$-x^2 \pa_x = \partial_r$ and $x \partial_{y^j}$ respectively\footnote{More
informally, a symbol such as $a(x,y,\nu,\mu)$ corresponds in the usual
polar co-ordinates $r,y$ to the pseudo-differential operator
$a(\frac{1}{r}, y, \frac{\partial}{\partial r}, \nabb)$.  The advantage of
these co-ordinates is that $\nu, \mu$ remain bounded for unit speed
geodesic flows.}; thus one can think of $\mu$ as an element of
$T^*_y(\partial M)$, the canonical one-form is $-{\nu dx}/{x^2} +
{\mu^j dy_j}/{x}$.  Following Melrose \cite{melrose}, we define the
\emph{scattering cotangent bundle} over the compact manifold $\overline{M}$
as the bundle $\Tscstar \Mbar$ whose sections are locally spanned over
$\mathcal{C}^\infty(\overline{M})$ by $dx/x^2$ and $dy/x$ (hence can be
paired with the \emph{scattering vector fields}, spanned by $-x^2 \pa_x$
and $x \partial_y$). Because of the form of the metric, the symbol
$\sigma(H)$ of the Schr\"odinger operator $H$ now has the form
$$ \sigma(H) = \frac{1}{2} (\nu^2 + h_{ij}(z) \mu^i \mu^j) = \frac{1}{2} (\nu^2 + |\mu|_{h(z)}^2)$$
and thus in the unit speed case $\sigma(H) = \frac{1}{2}$, $\nu$ and $\mu$ will have 
combined magnitude equal to 1.  

The geodesic flow on $\Tscstar \Mbar$ can then be written explicitly as
\begin{equation}\begin{aligned}
\frac{d}{ds} x &= -x^2 \nu \\ 
\frac{d}{ds} y^j &= x h^{jk} \mu_{k} \\ 
\frac{d}{ds} \nu &= x h^{jk}
\mu_j \mu_k + \frac{1}{2} x^2 \frac{\partial h^{jk}}{\partial x} \mu_j
\mu_k \\ 
\frac{d}{ds} \mu_j &= -x \mu_j \nu - \frac{1}{2} x \frac{\partial
h^{kl}}{\partial y^j} \mu_k \mu_l;
\label{ham-flow}\end{aligned}\end{equation}
recall that $h^{jk}(x,y)$ is, for small $x,$ a family of metrics in $y$
with $x$ a smooth parameter; it is only in the perfectly conic case that
there is no $x$-dependence.

This geodesic flow of course still preserves $\sigma(H)$, and is still homogeneous with respect to the scaling
$(x(s),y(s),\nu(s),\mu(s)) \mapsto (x(\lambda s), y(\lambda s), \lambda \nu(\lambda s), \lambda \mu(\lambda s))$.
This geodesic flow equation can be used to define the distance $d_M((x',y'),(x'',y''))$ provided $(x',y')$ and
$(x'',y'')$ are sufficiently close; we will quantify what `sufficiently close' means later.  

To gain some intuition as to how the distance $d_M$ behaves on $\Omega(y_0,\eta,r_0)$, let us
first consider the case $h^{jk}(x,y) = h^{jk}(y)$ when $M$ is perfectly conic near infinity; to distinguish this 
from the asymptotically conic situation we shall write $d_\conic = d_{\conic,h}$ instead of $d_M$, and write
the metric $g$ as $g_{\conic, h}$.  In this case 
the equations for geodesic flow
are equivalent to those on a plane, for the following geometric reason.  Consider a geodesic
$\gamma$ in $B_{\partial M}(y_0,\eta)$, and consider the truncated cone 
$\{ (x,y): 0 < x < \epsilon_0, y \in \gamma \}$ over that geodesic.
In the perfectly conic case, this truncated cone is both perfectly flat and totally 
geodesic\footnote{Indeed, one can think of $\gamma$ as being isometric to an arc in $S^1$, so the truncated cone is then isometric to a truncated sector in $\R^2$.}, and in particular is isometric to a subset of the plane.  Conversely, 
it is easy to show that 
there is a unique minimal geodesic connecting $(x',y')$ and $(x'',y'')$ which lies in the
truncated cone over the geodesic in $\partial M$ connecting $y'$ and $y''$.  From plane geometry we thus see that
we have the \emph{cosine rule}\footnote{This rule can also be deduced from \eqref{ham-flow}; the point is that
the $\frac{\partial h^{ij}}{\partial x}$ term in the equation for $\frac{d}{ds} \nu$ disappears in the perfectly
conic case, and so the geodesic flow equation decouples into a (rescaled) geodesic flow equation on $\partial M$
for $y, \mu$, and an explicitly solvable evolution for $x$, $\nu$.  We omit the details.}
\begin{equation} d_\conic(z',z'')^2 = (r')^2 + (r'')^2 - 2 r' r'' \cos d_{\partial M}(y',y'')
\label{exact-conic-dist}\end{equation}
for $z', z'' \in \Omega(y_0,\eta,r_0)$.  
Note in particular that we thus expect symbol-type regularity estimates on $d_\conic$, in the sense
that each derivative in the $z'$ and $z''$ variable gains a power of $d$; this is of course what happens in the
Euclidean case, and we will be able to show that it also happens in asymptotically conic manifolds near infinity.
This is quite plausible when $r'$ and $r''$ are of comparable magnitude; the most delicate issue will be proving
these sorts of bounds when $r'$ and $r''$ are both large but far apart from each other, which seems to be a geometric
necessity in order to not lose any derivatives in Theorem \ref{main}.

This completes our discussion of the perfectly conic case.  We now consider
the distance function $d_M$
of an asymptotically conic manifold on $\Omega(y_0, \eta, r_0)$.  We shall always take $\eta$ to be sufficiently small
and $r_0$ to be sufficiently large. 

We first show that globally geodesics between points in $\Omega(y_0,\eta,r_0)$ stay in the scattering region.

\begin{lemma}\label{min-exists} If $\eta$ is sufficiently small and $r_0$ sufficiently large, then for every pair of points $z', z''$ in $\Omega(y_0, \eta, r_0)$, every globally minimizing geodesic between them lies entirely inside $\Omega(y_0, 5\eta, r_0/2)$.
\end{lemma}

\begin{proof} Without loss of generality we may take $r'' \geq r'$.
Suppose for contradiction that the globally minimizing curve connecting $z'$ and $z''$ contains a point $(r,y)$ in the complement of $\Omega(5\eta, r_0/2)$. Then either $r < r_0/2$ or $d_{\partial M}(y, y_0) > 5\eta$. 

First suppose that $r < r_0/2$.   Using the crude estimate $g \geq dr^2$,
we thus obtain the lower bound
\begin{equation} d_M(z',z'') \geq (r'-r_0/2)+ (r''-r_0/2) \geq r''.
\label{lb1}\end{equation}
On the other hand, if $r_0$ is large enough then we have the pointwise
bound $h_{ij}(1/r,y)$ $\leq 4 h_{ij}(0,y)$ for all $r > r_0/10, y \in
\partial M$, and hence $g \leq g_{\conic, 4h}$ in ths region.  Since the
minimal geodesic connecting $z'$ to $z''$ in the $g_{\conic,4h}$ metric
lies in the region $r > r_0/10$, we thus obtain from \eqref{exact-conic-dist}
\begin{equation}
d_M(z', z'') \leq d_{\conic,4h}(z',z'') \leq \sqrt{(r')^2 + (r'')^2 - 2r'r''\cos 4\eta},
\label{over}\end{equation}
which contradicts \eqref{lb1} if $\eta$ is sufficiently small.

Next suppose that $d_{\partial M}(y, y_0) > 5\eta$. Then the distance from
$y$ to both $y'$ and $y''$ must be at least $4\eta$. Using the comparison
$h_{ij}(0,y)/4 \leq h_{ij}(1/r,y) \leq 4h_{ij}(0,y)$, which holds for $r >
r_0/2$ if $r_0$ is sufficiently large, we have the pointwise bound $dr^2 +
h_{ij}(0,y)/4$ $\leq g \leq dr^2 + 4h_{ij}(0,y)$. Comparison with this
smaller conic metric shows that the length of the curve must be at least
\begin{equation}
\sqrt{(r')^2 + r^2 - 2r'r\cos 2\eta} + \sqrt{(r^2 + (r'')^2 - 2rr''\cos 2\eta},
\label{under}\end{equation}
since the broken geodesic in the small conic metric that goes between these
points has at least this length. For $\eta < \pi/12$ this is larger than
\eqref{over}, as can be seen by consideration of the lengths of a plane
Euclidean triangle with vertices having polar coordinates $(r',0)$, $(r'',
4\eta)$ and $(r', 2\eta)$. This is a contradiction.
\end{proof}

We now use a symplectic argument to deduce regularity of
the distance function $d_M(z', z'')$ where $(z', z'')$ lie in some cone $\Omega(y, \eta, r_0)$ for fixed $(\eta, r_0)$.  To do this, we fully exploit the smoothness of the metric function $h_{ij}(x,y)$ in scattering coordinates $(x,y)$. We shall need to work on the a blown up version of the compactified double space $\overline{M}^2$ in order to do this. 
 
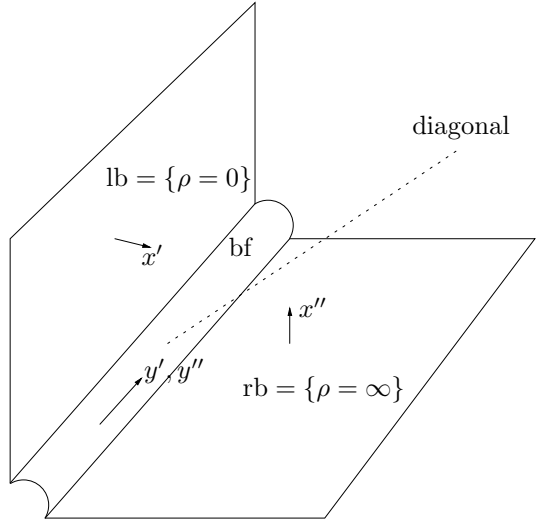
\begin{figure}
\setlength{\unitlength}{0.0003in}
\begingroup\makeatletter\ifx\SetFigFont\undefined%
\gdef\SetFigFont#1#2#3#4#5{%
  \reset@font\fontsize{#1}{#2pt}%
  \fontfamily{#3}\fontseries{#4}\fontshape{#5}%
  \selectfont}%
\fi\endgroup%
{\renewcommand{\dashlinestretch}{30}
\begin{picture}(9024,8889)(0,-10)
\put(237.000,237.000){\arc{874.643}{4.1720}{6.8236}}
\put(4437.000,5037.000){\arc{874.643}{4.1720}{6.8236}}
\path(12,612)(4212,5412)
\path(612,12)(4812,4812)
\path(612,12)(5412,12)(9012,4812)(4812,4812)
\path(12,612)(12,4812)(4212,8862)(4212,5412)
\dashline{60.000}(2712,3012)(7662,6312)
\path(4812,3012)(4812,3612)
\blacken\path(4842.000,3492.000)(4812.000,3612.000)(4782.000,3492.000)(4842.000,3492.000)
\path(1812,4812)(2412,4662)
\blacken\path(2288.307,4662.000)(2412.000,4662.000)(2302.859,4720.209)(2288.307,4662.000)
\path(1554,1625)(2262,2412)
\blacken\path(2204.046,2302.724)(2262.000,2412.000)(2159.440,2342.852)(2204.046,2302.724)
\put(6912,6612){diagonal}
\put(1662,5712){$\lb = \{ \rho = 0 \}$}
\put(4002,2112){$\rb = \{ \rho = \infty \}$} 
\put(3762,4512){$\blf$}
\put(4962,3462){$x''$}
\put(2262,4362){$x'$}
\put(2312,2412){$y',y''$}
\end{picture}
}
\caption{The blown up space $\overline{M}^2_b$}\label{fig:m2b}
\end{figure}

Let $\overline{M}^2_b$ denote the blown-up compactified manifold consisting
of $\overline{M}^2$ with the boundary $(\partial M)^2$ blown up in the
sense of \cite{melrose} (shown in Figure~\ref{fig:m2b}).  We define $\rho =
x'/x''$ and label the hypersurfaces of $\overline{M}^2_b$ by $\lb$ (`left
boundary'), $\rb$ (`right boundary') and $\blf$ (`blown-up face') according
as they arise from the faces $x' = 0$, $x'' = 0$ or $x' = x'' = 0$ of
$\overline{M}^2$, respectively. Let $U$ be any neighbourhood in
$\overline{M}^2_b$ of the closure of the diagonal  $\{ (z,z): z
\in M \} \cup \{ w \in \blf: \rho(w) = 1, y'(w) = y''(w) \}$.  For
instance, $U$ could be the region $\{ (z',z'') \in M: d_M(z',z'') <
c(\langle z' \rangle + \langle z'' \rangle) \}$ for some small $0 < c \ll
1$, together with the boundary of this region on the blown up face $\blf$.

\begin{proposition}
Let $\Omega(\eta, r_0) \subset \overline{M}^2_b$ be defined by
$$
\Omega(\eta, r_0) = \{ x',x''< r_0^{-1},\ d_{\partial M}(y',y'')<\eta \}
$$
Then for $r_0$ sufficiently large and $\eta$ sufficiently small,
\begin{equation}
\text{ the function } \frac{d_M(z',z'')}{r' + r''} \text{ is in } \CI(\Omega(\eta, r_0) \setminus U).
\label{wtp}\end{equation}
\end{proposition}

\begin{remark} Notice that $U$ includes a \emph{conic} neighbourhood of the diagonal near infinity, that is, a set of the form 
\begin{equation}
\{ (r',y',r'',y'') \mid \big| \frac{r''}{r'} - 1 \big| + d_{\partial M}(y',y'') < \epsilon \}
\label{diag-conic-nbhd}\end{equation}
for some $\epsilon > 0$, because of the topology of $\overline{M}^2_b$.
\end{remark}

\begin{proof}
To prove the claim, we use a symplectic technique; namely, we use the
function $\Phi = d_M(z',z'')$ as a generating function for a Lagrangian
submanifold. This function satisfies the eikonal equation
\begin{equation}\label{eikonal}
|\nabla_{z'} \Phi|^2 = 1,\end{equation} and vanishes at the diagonal.
  Consequently the graph of $d\Phi$ will be tangent to the flow in $\Tscstar(\Mbar) \times \Tscstar(\Mbar)$, lifted to a vector bundle over $\Mbar^2_b$, generated by the energy $\frac{1}{2}(|\zeta'|^2_g - 1)$. This flow is given by \eqref{ham-flow} for the singly-primed variables, while $z''$ and $\zeta''$ are stationary under the flow. Using coordinates $(\rho, y', x'', y''; \nu', \mu', \nu'', \mu'')$, the vector field and the evolution of $\Phi$ is given by 
\begin{equation}\begin{gathered}\begin{aligned}
\frac{d}{ds} \rho &= -x' \rho \nu \\ 
\frac{d}{ds} {(y')}^i &= x' h^{ij} \mu'_j \\ 
\frac{d}{ds} \nu' &= x' h^{ij} \mu'_i \mu'_j + \frac{(x')^2}{2} \frac{\partial h^{ij}}{\partial x'} \mu'_i \mu'_j \\ 
\frac{d}{ds} \mu'_k &= -x' \mu'_k \nu' - \frac{x'}{2}  \frac{\partial h^{ij}}{\partial (y')^k} \mu'_i \mu'_j.
\end{aligned}\qquad \begin{aligned}
\frac{d}{ds} x'' &= 0 \\
\frac{d}{ds} y'' &= 0 \\ 
\frac{d}{ds} \nu'' &= 0 \\
\frac{d}{ds} \mu'' &= 0 
\end{aligned} \qquad
\frac{d}{ds} \Phi = 1. \end{gathered}
\label{hvf}\end{equation}
Here $h_{ij}$ is always a function of the $(x',y')$ variables. 
Geometrically we consider the submanifold $L$ determined by the `initial condition' that 
\begin{equation}
\{ x' = x'', \ y' = y'', \ \mu' = - \mu'', \ \nu' = -\nu'', \ h^{ij}(x',y')\mu'_i \mu'_j + \nu^2 = 1 \} \subset L, 
\label{init-cond}\end{equation}
and $L$ is invariant under the forward flow of \eqref{hvf} (i.e.\ for
positive `time'). This submanifold is a smooth Lagrangian manifold, with
boundary, given by the graph of $d\Phi$, at least near the diagonal. (It is
perhaps counterintuitive that $L$ is smooth near the diagonal, since $\Phi$
is singular there. The singularity of $\Phi$ is reflected in the fact that
$L$ projects diffeomorphically to $\overline{M}^2_b$ near the diagonal, but
not at the diagonal. The inverse image of a point on the diagonal is a
$(n-1)$-sphere consisting of all possible unit directions at that
point. Smoothness of $L$ is reflected in the fact that $d\Phi$ does not
blow up at the diagonal, rather it becomes `multivalued'.)

Smoothness of $L$ persists as long as the vector field of which it is the flowout remains smooth and nonvanishing. We want to show that the function $\Phi/(r' + r'')$ extends to a smooth function on $\Omega_b$ (including its boundary at infinity). To do this, we shall multiply the vector field \eqref{hvf} by a factor, which does not change the flowout. The new vector field will be smooth \emph{and nonvanishing up to the boundary} in suitable coordinates, and this will imply that $L$ is a smooth submanifold in these coordinates. 

By symmetry of $\Phi$ under interchange of $z'$ and $z''$, it is sufficient to work in the region $\{ r'' < 2r' \} = \{ \rho < 2 \}$. We shall refer to $\Omega(\eta, r_0) \cap \{ \rho < 2 \}$ as the `region of interest'. In the region of interest it is equivalent to require that $\kappa = x' \Phi = \Phi/r'$, as opposed to $\Phi/(r' + r'')$, is smooth on $\Omega_b$. Its evolution is given by 
\begin{equation*}
\frac{d}{ds} \kappa = x' \frac{d}{ds} \Phi + \phi \frac{d}{ds} x' = x' - x' \kappa \nu'.
\end{equation*}

As a first step we multiply the flow \eqref{hvf} by $(x')^{-1}$, obtaining the new flow
\begin{equation}\begin{gathered}\begin{aligned}
\frac{d}{ds'} \rho &= -\rho \nu \\ 
\frac{d}{ds'} {(y')}^i &= h^{ij} \mu'_j \\ 
\frac{d}{ds'} \nu' &= h^{ij} \mu'_i \mu'_j + \frac{x'}{2} \frac{\partial
  h^{ij}}{\partial x'} \mu'_i \mu'_j \end{aligned}
\qquad \qquad
\begin{aligned}
\frac{d}{ds'} \mu'_k &= - \mu'_k \nu' - \frac{1}{2}  \frac{\partial
  h^{ij}}{\partial (y')^k} \mu'_i \mu'_j \\
\frac{d}{ds'} \kappa &= 1 - \kappa \nu', \end{aligned}\end{gathered}
\label{hvf2}\end{equation}
where $ds' := x' ds$.  This flow is defined for all `time', but the vector
field vanishes where $\rho = 0, \mu = 0$, so we need to desingularize about
this set (i.e.\ blow it up). Let us introduce the coordinates $M_i =
\mu'_i/\rho$, $K = (\kappa - 1)/\rho$, and consider the flow just on the
energy surface $\{ (\nu')^2 + h^{ij} \mu'_i \mu'_j = 1 \}$. (Introduction
of $M$ is equivalent to working on the bundle $\Tscstar(\Mbar) \times
\Tscstar(\Mbar)$ lifted to $\Mbar^2_b$ and then blown up at $\{ \rho = 0,
\mu' = 0 \}$).  In terms of these quantities, we get, after dropping the
primes on $\nu$ and $y$ and changing to a new parameter given by $ds''=\rho
\, ds'$
\begin{equation}
\begin{gathered}
\begin{aligned}
\frac{d}{ds''} \rho &= - \nu \\
\frac{d}{ds''} {y}^i &= h^{ij} M_j \\ 
\frac{d}{ds''} \nu &= \rho h^{ij} M_i M_j + \frac{\rho^2 x''}{2}
\frac{\partial h^{ij}}{\partial x} M_i M_j  \end{aligned}
\qquad \qquad
\begin{aligned}
\frac{d}{ds''} M_k &=- \frac{1}{2}  \frac{\partial h^{ij}}{\partial y^k} M_i M_j\\
\frac{d}{ds''} K &= \frac{h^{ij}M_i M_j}{1 + \nu}.\end{aligned}\end{gathered}
\label{hvf4}\end{equation}
This remains nonsingular as long as $|M|$ remains bounded and $1 + \nu \neq
0$. For the remainder of this paragraph, let us consider what happens at
the face bf = $\{ x'' = 0 \}$, to which \eqref{hvf4} is tangent. At bf, the
flow is identical to the flow for the perfectly conic metric $dr^2 + r^2
h(0)$. Since $d_{s''} (h^{ij}M_i M_j) = 0$ and initially, $-1 \leq |M|_h \leq 1$,
$M$ certainly remains bounded. As for $\nu$, we have $\nu' \geq 0$, and
under the flow \eqref{hvf4}, $d_{s''} \rho = -\nu$.  Thus if $\rho$
becomes smaller than $1$ along a trajectory, then at some point $\nu \geq
0$ and therefore $\nu > 0$ thereafter along the trajectory. Hence at $\rho
= 0$, $\nu \geq 0$ and there is no singularity in the equation for
$d_{s''} K$.  Hence the flow at bf is well behaved in the region of
interest, and it is not hard to show that it exits the region of interest
in finite `time' (measured by the parameter $s''$).

By the theory of ODEs depending on a smooth parameter, for sufficiently small initial values of $x''$ and for initial values given by \eqref{init-cond}, $|M|$ remains bounded and $1 + \nu$ strictly positive in the region of interest, and the region of interest is exited in finite `time'. We see then that for small $x''$, the submanifold $L$ is smooth in the region of interest, and $K$, and hence also $\kappa$, is a smooth function on $L$. 

Next we claim that in $\Omega(\eta, r_0) \cap \{ \rho < 2 \}$, for $\eta$ sufficiently small and $r_0$ sufficiently large, $(y', y'', \rho, x'')$ furnish coordinates on $L$. At $L \cap \{ x'' = 0 \}$, $\Phi$ is given by the explicit conic formula \eqref{exact-conic-dist} when $y'$ and $y''$ are close:
\begin{equation*}
\Phi = \big(  (r')^2 + (r'')^2 - 2r' r'' \cos d_{\partial M}(y',y'') \big)^{1/2} = r' \big( 1 + \rho^2 + 2 \rho  \cos d_{\partial M}(y',y'') \big)^{1/2}.
\end{equation*}
This is smooth away from 
\begin{equation}
\rho = 1, \quad y' = y''.
\label{sing}\end{equation}
Therefore the graph of $d\Phi$, restricted to $x'' = 0$, has $(y', y'',
\rho)$ as smooth coordinates on any open set excluding \eqref{sing}, and
hence, their differentials are linearly independent on $L \cap \{ x'' = 0
\}$. Also, the differential of $x''$ is nonzero when restricted to $L$ at
$x'' = 0$. This follows from the fact that the differential is zero at the
initial hypersurface \eqref{init-cond}, and that the Lie derivative of
$dx''$ with respect to \eqref{hvf4} vanishes. Therefore, $L$ is
\emph{projectable} --- in other words, $dy', dy'', d\rho, dx''$ have
linearly independent differentials when restricted to $L$ at bf, on any
open set excluding \eqref{sing}, and therefore form coordinates on $L$ ---
\emph{uniformly up to the corner } $\{ x'' = \rho = 0 \}$. By continuity,
this remains true on a neighbourhood of $L \cap \{ x'' = 0 \}$, which
covers an $\Omega(\eta, r_0)$ for sufficiently small $\eta$ and large
$r_0$. Therefore, $\kappa$ is a smooth function of $(y',y'', x'', \rho)$ on
$\Omega(\eta, r_0) \setminus U$ for sufficiently small $\eta$ and large
$r_0$. Hence for $r_0$ sufficiently large and $\eta$ sufficiently small,
there is exactly one geodesic
between any two points $(r',y''), (r'',y'')$ with $r',r'' > r_0$ and
$d_{\partial M}(y',y'') < \eta$ that lies wholly inside $\Omega(5\eta,
r_0/2)$. Lemma~\ref{min-exists} shows that every globally minimizing
geodesic lies in this region. Therefore, all geodesics lying wholly within
$\Omega(5\eta, r_0/2)$ are globally minimizing geodesics. This proves that
$\Phi = r' \kappa$ is the genuine distance function on $M$. This completes
the proof of the proposition.
\end{proof}

We are now able to approximate the asymptotically conic metric $d_M$ rather
precisely by the perfectly conic analogue $d_{\conic} = d_{\conic,h}.$

\begin{proposition}\label{cor:dist}  We have 
\begin{equation}
d_M(z',z'') = d_{\conic}(z',z'') + e(z',z''),
\label{dist-approx}\end{equation}
where $d_{\conic}$ is the distance function for the perfectly conic metric
$dr^2 + r^2 h_{ij}(0)dy^i dy^j$ (which is then given by \eqref{exact-conic-dist}), and the error $e(z',z'')$ is $C^\infty$ on $\Omega(\eta, r_0) \setminus U$. 
\end{proposition}

\begin{proof} We claim that $\kappa$ and $\kappa_{\conic}$ agree at each boundary hypersurface of $\Mbar^2_b$. In fact, at lb $ = \{ \rho = 0 \}$ they are both equal to $1$, while at bf, they solve the same ODE with the same initial conditions. Equation \eqref{dist-approx} follows immediately. 
\end{proof}

We can now prove Lemma \ref{symbol-est}.  We first prove the claim for
perfectly conic metrics outside of $U$. 

\begin{lemma}\label{perfect-symbol}  Let $U$ be the region
$U := \{ (z',z'') \in M: d_M(z',z'') < (\langle z' \rangle + \langle z'' \rangle)/100 \}$.  Then we have
$$ |\nabla_{z'}^{m'} \nabla_{z''}^{m''} \nabla_{z',z''} d_\conic(z',z'')| \leq C_{m',m''} \langle z' \rangle^{-m'}
\langle z'' \rangle^{-m''}$$
for all $m',m'' \geq 0$ and $z',z'' \in \Omega(\eta,r_0) \backslash U$.
\end{lemma}

\begin{proof}
We of course use \eqref{exact-conic-dist}.  Without loss of generality we
may take $\langle z' \rangle \geq \langle z'' \rangle$; since we are
outside $U$, we thus have $d_\conic$ comparable to $\langle z' \rangle$.
We begin by computing derivatives of $d_\conic^2$.  Observe that $\cos
d_{\partial M}(y',y'')$ is a smooth function of $d_{\partial M}(y',y'')^2$.
This function is in turn smooth for $y', y'' \in B_{\partial M}(y_0,\eta)$
for $\eta$ sufficiently small (see e.g.\ \cite{jost}), so one can then
verify the bounds
$$|\nabla_{z'}^{m'} \nabla_{z''}^{m''} d_\conic^2(z',z'')| \leq C_{m',m''} \langle z' \rangle^{1-m'}
\langle z'' \rangle^{1-m''}$$
when $m' \geq 0$ and $m'' \geq 1$ (because the $(r')^2$ term vanishes); for pure $z'$ derivatives we have the 
slightly different formula
$$ |\nabla_{z'}^{m'} d_\conic^2(z',z'')| \leq C_{m'} \langle z' \rangle^{2-m'}.$$
Expanding the latter inequality using the Leibnitz rule, one can then prove
$$ |\nabla_{z'}^{m'} d_\conic(z',z'')| \leq C_{m'} \langle z' \rangle^{1-m'}$$
inductively for all $m' \geq 0$, and then by expanding the former inequality using the Leibnitz rule one can also prove
$$|\nabla_{z'}^{m'} \nabla_{z''}^{m''} d_\conic(z',z'')| \leq C_{m',m''} \langle z' \rangle^{-m'}
\langle z'' \rangle^{1-m''}$$
inductively for all $m' \geq 0$ and $m'' \geq 1$.  The claim follows.
\end{proof}

Next, we verify symbol estimates on a compact portion of the manifold.

\begin{lemma}\label{compact-symbol}  Let $z_0 \in M$, and let $\eta > 0$ be sufficiently small depending on $z_0$.
Then we have
$$ |\nabla_{z',z''}^m d_M(z',z'')| \leq C_{z_0,\eta,m} d_M(z',z'')^{1-m}$$
for all $m \geq 0$ and $z', z'' \in B(z_0,\eta)$.
\end{lemma}

\begin{proof}  We work in normal co-ordinates $z^j$ around $z'$, so $z^j(z') = 0$ and $\sum_j |z^j(z'')|$ is comparable
to $d_M(z',z'')$.  Write $\tau := d_M(z',z'')/\eta$, and rescale the
co-ordinates by $\tau$ (cf.\ Lemma \ref{lemma:lapsquareddist}), giving rise
to a new manifold $M_\tau$ with metric $g^{jk}_\tau(z) = g^{jk}(\tau z)$.
In the ball $\{ \sum_j |z^j| < C\eta \}$, these metrics $g_\tau$ vary
smoothly as $\tau$ varies over the compact set $[0,1]$, and the distance
functions $d_{M_\tau}(z',z'')$ also vary smoothly in the region where
$\sum_j |z^j(z')| < C\eta$, $\sum_j |z^j(z'')| < C\eta$, and $\sum_j
|z^j(z') - z^j(z'')| > c\eta$.  In particular we have estimates of the form
$$ |\nabla_{z',z''}^m d_{M_\tau}(z',z'')| \leq C_{z_0,\eta,m}$$
uniformly in $\tau$ in this region.  The claim then follows by undoing the scaling.
\end{proof}

Now we combine these two estimates to obtain Lemma \ref{symbol-est}.

\begin{proof}[Proof of Lemma \ref{symbol-est}]  We let $U$ be the region 
$\{ (z',z'') \in \Omega(\eta,r_0): d_M(z',z'') < c(r'+r'') \}$ for some small $c > 0$, together with the boundary of this region on $\blf$.
Outside of $U$, we have $d_M(z',z'')$ comparable to $r'+r''$, and
the symbol bound on $d_M$ follows from \eqref{dist-approx}, since the error $e$ being smooth on the compactified manifold automatically obeys symbol estimates (indeed it obeys estimates which are one order of $r'+r''$ better
than required), while the symbol estimates $d_\conic$ follow from Lemma \ref{perfect-symbol}.  Inside the
region $U$, we can then rescale as in Lemma \ref{lemma:lapsquareddist} to rescale $z'$ and $z''$ to a fixed compact set,
and then apply Lemma \ref{compact-symbol}, if we choose $c$ sufficiently small.  Note that the claim concerning
the nature of the singularity of $d_M$ at the diagonal $z'=z''$ just follows from the fact that $d_M^2$ is smooth
and vanishes to second order at the diagonal (see e.g.\ \cite{jost}).
\end{proof}

It remains to prove Proposition \ref{full-ang-est}.  We first show this in the case of perfectly conic
manifolds (in which the non-angular error terms do not appear).

\begin{lemma}\label{conic-est} For the perfectly conic distance $d_{\conic}(z',z'')$, we have the estimate\footnote{We caution the reader that there are \emph{two} `time' variables in play here; the variable $t$ used to parameterize
the moving particles $z'(t)$ and $z''(t)$, and also the implicit variable $s$ used to parameterize the geodesic
connecting $z'(t)$ and $z''(t)$.} 
\begin{equation}
\frac{d^2}{dt^2} d_{\conic}(z'(t), z''(t)) \geq - C\Big( \frac{|v'_{\perp}|^2}{r'} + \frac{|v''_{\perp}|^2}{r''} \Big),
\label{conic-ang-est}\end{equation}
for $d_{\partial M}(y',y'')$ less than the injectivity radius of $(\partial
M, h(0))$, where $z'(t)$ moves along a geodesic (for the conic metric) with
initial condition $(z', v')$, and $z''(t)$ moves along a geodesic with
initial condition $(z'', v'')$. Here $v'_{\perp}$ resp.\ $v''_{\perp}$
denotes the component of $v'$ resp.\ $v''$ perpendicular to both
$\partial_r$ and $\frac{d}{ds}\gamma$, the tangent vector to the geodesic
$\gamma$ from $z'$ to $z''$.
\end{lemma}
Note that, in the notation of Proposition~\ref{full-ang-est},
$\abs{v'_\perp} \leq \abs{v'_\ang}$, $\abs{v''_\perp} \leq \abs{v''_\ang}$,
hence this result implies Proposition~\ref{full-ang-est} in the perfectly
conic case.

\begin{proof}
Let us decompose the tangent vectors $v' = v'_{\perp} + v'_{\para}$, $v' =
v''_{\perp} + v''_{\para}$, where $v'_{\para}$ is in the span\footnote{If
$z'$ and $z''$ lie on the same radial arc, then $\gamma'$ and $\partial_r$
are parallel, but the argument still works (e.g.\ by adjoining an arbitrary
vector to the span of $\partial_r$ to define the $v'_{\para}$ component) in
this case.} of $\gamma'$ and $\partial_r$, and $v'_{\perp}$ is
perpendicular to the span of $\gamma'$ and $\partial_r$, and similarly for
$v''$.  Clearly the left-hand side of \eqref{conic-ang-est} is a bilinear
form in these velocity variables.

First suppose that $v'_{\perp} = v''_{\perp} = 0$. The geodesic $\gamma(s)$
lies in a flat plane which is totally geodesic. Hence, in this case all
geodesics $\gamma(t)$ from $z'(t)$ to $z''(t)$ lie in a flat plane, and the
truth of \eqref{conic-ang-est} follows from its truth in Euclidean $\R^2$, which was shown in Section~\ref{morawetz-sec}. 

Next suppose that $v'_{\para} = v''_{\para} = 0$. In this case, we use the formula \eqref{exact-conic-dist} and compute explicitly, exploiting the fact that in this case 
$$\frac{d}{dt} d_M(z', z'')|_{t=0} = \frac{d}{dt} r'(t)|_{t=0} = \frac{d}{dt} r''(t)|_{t=0} = 
\frac{d}{dt} d_{\partial M}(y',y'')|_{t=0} = 0.$$
For a perfectly conic metric we have from \eqref{ham-flow} that
$$\frac{d^2}{dt^2} (r')^2(t) = 2|\mu'|^2, \quad \frac{d^2}{dt^2} r'(t) = |\mu'|^2/r',$$
hence $d^2/dt^2 d_{\conic}(z',z'')$ is given by 
\begin{equation}\begin{gathered}
\frac1{d_{\conic}(z',z'')} \Bigg( 
\frac{d^2}{dt^2} ((r')^2) + \frac{d^2}{dt^2} ((r'')^2)- 2 \big( (\frac{d^2}{dt^2} r') r'' + (\frac{d^2}{dt^2} r')r'' \big) \cos d_{\partial M}(y',y'') \\ + 2 r' r''  \frac{d^2}{dt^2} {\cos d_{\partial M}}(y',y'') \Bigg)
\end{gathered}\label{est-1}\end{equation}
Note that $(\cos d_{\partial M}(y',y'')) =  (1 -  d_{\partial M}(y',y'')^2/2 + O(d_{\partial M}^4) )$ is a smooth function of $(y', y'')$. A computation in normal coordinates around (say) $y'$ shows that 
\begin{equation*}
\frac{d^2}{dt^2} (\cos d_{\partial M}(y',y'')) \geq -  C\Big( \frac{|\mu'|}{r'} + \frac{|\mu''|}{r''} \Big)^2 d_{\partial M}(y',y''). 
\end{equation*}
Since $d_{\partial M}(y',y'')$ is bounded above, we can estimate \eqref{est-1} from below by
\begin{equation}
\frac1{d_{\conic}(z',z'')} \Bigg( 2 (|\mu'|^2 +  |\mu''|^2) - 2 \big( \frac{r''|\mu'|^2}{r'} + \frac{r'|\mu''|^2}{r''} \big) - Cr'r'' \big( \frac{|\mu'|}{r'} + \frac{|\mu''|}{r''} \big)^2  \Bigg)
\label{est-2}\end{equation}
We rewrite this as 
\begin{equation}
\frac1{d_{\conic}(z',z'')} \Bigg( 2 \big(\frac{(r'-r'')|\mu'|^2}{r'} +  \frac{(r''-r')|\mu''|^2}{r''}\big)  - C\big( \frac{r''|\mu'|^2}{r'} + \frac{r'|\mu''|^2}{r''} \big)  \Bigg)
\label{est-3}\end{equation}
Since $|r'-r''| < d_{\conic}(z',z'')$, the first term satisfies the
estimate \eqref{conic-ang-est}. For the second term, if $r' \leq 2r'' \leq
4r'$ then it automatically satisfies the estimate \eqref{conic-ang-est},
while if the ratio of $r'$ and $r''$ is large then $d_{\conic}(z',z'')$ is
comparable to the \emph{larger} of $r'$ and $r''$, showing that in this
case also the second term satisfies the estimate \eqref{conic-ang-est}.

Finally consider the case when one vector, say $v'$, is angular, $v'_{\para} = 0$, and the other is parallel, $v''_{\perp} = 0$. In this case, 
we write 
\begin{equation*}
\frac{d^2}{dt^2} d_{\conic} = \nabla^2_{(v',v'')} d_{\conic} = \nabla^2_{(v',0)} d_{\conic} + 2 \nabla_{(v',0)}\nabla_{(0,v'')} d_{\conic} + \nabla^2_{(0,v'')} d_{\conic}.
\label{cross}\end{equation*}
and note that the geodesic $\gamma_t$ from $z'(0)$ to $z''(t)$ lies in the
flat plane containing $\gamma$ for all $t$. Hence $v'=v'_{\perp}$ is
perpendicular to $\gamma_t$ for all $t,$ and we obtain
$$
\nabla_{v'} d(z', z''(t))=0 \quad \text{for all } t.
$$ Differentiating with respect to $t$ then shows that the cross term in
\eqref{cross} vanishes. The estimate \eqref{conic-ang-est} then follows,
since the terms $\nabla^2_{(v',0)} d_{\conic}$ and $\nabla^2_{(0,v'')}
d_{\conic}$ have already been treated.
\end{proof}

\begin{lemma}\label{lemma:ang} Suppose that $((r',y'), (r'',y'')$ are in a small neighbourhood $U$ of the diagonal  in $\overline{M}^2_b$, say in the set \eqref{diag-conic-nbhd} for $\epsilon > 0$ sufficiently small. Then the
estimate \eqref{ang-est2} holds.
\end{lemma}

\begin{proof}
We begin by scaling the problem so that the initial point $z'(0) = (r'(0),
y'(0))$ is at $r = 1$. Let $\tau^{-1} = r'(0)$. We scale by factor $\tau$;
that is, we consider the metric $g_\tau = dr^2 + r^2h(\tau/r, y)$ and
replace the geodesic $(r(t), y(t))$ by $(\tau r(t), y(t))$ which is a
(reparametrized) geodesic with respect to $g_{\tau}$. An equivalent
formulation of the lemma is to prove that
\begin{equation*}\begin{gathered}
\frac{d^2}{dt^2} d_{\conic}(z'(t), z''(t)) \geq - C\Big( |v'_{\perp}|^2 +
|v''_{\perp}|^2 \Big) \\ - C \tau \Big( |v'|^2 + |v''|^2 \Big)
\end{gathered}\end{equation*}
for all geodesics joining points $(r',y')$, $(r'',y'')$ satisfying $r' = 1$, $|r'' - 1| < \epsilon$, $d_{\partial M}(y',y'') < \eta$. 

To do this, we return to the second variation formula, and use the fact that
$$
\frac{d^2}{dt^2} d_M(z'(t), z''(t)) \geq  -\int_0^{d(z'(0),z''(0))} \langle R(T,J)T, J \rangle ds
$$ where $J$ is the Jacobi field corresponding to the family of unit speed
geodesics $\gamma_t(s)$ between $z'(t)$ and $z''(t)$ and $T=d\gamma/ds.$
Let the boundary value of $J$ at $z'$ be decomposed as above in to
$J=J'_{\perp}+J'_{\para}$, and define $J''_{\perp}$ and $J''_{\para}$
similarly. Let $J_1$ be the Jacobi field with boundary values $J'_{\perp}$,
$J''_{\perp}$ and let $J_2$ be the Jacobi field with boundary values
$J'_{\para}$, $J''_{\para}$.  For the perfectly conic metric at $\tau = 0$,
$J_2$ is in the plane spanned by $\frac{d}{ds} \gamma$ and
$\partial_r$. Hence in the general case $J_2(\tau)$ is in the plane spanned
by $\frac{d}{ds} \gamma$ and $\partial_r$ up to $O(\tau)$. This means that
for any $V,$
\begin{equation*}
\langle R(T, J_2) T, V \rangle \text{ and } \langle R(T, V) T, J_2 \rangle = O(\tau |J_2| |V|).
\end{equation*}
Hence, using the fact that $d(z'(0), z''(0))$ is bounded on our rescaled
region of interest,
\begin{equation}
\frac{d^2}{dt^2} d(z'(t), z''(t)) \geq  -C(\norm{J_1}_\infty^2+\tau^2 \norm{J_2}_\infty^2).
\label{J-est}
\end{equation}
For $\tau \in [0,1]$, we can estimate the length of a Jacobi field
corresponding to $v=(v',v'')$ by
\begin{equation}
\| J \|_{\infty}^2 \leq C (|v'|^2 + |v''|^2)
\label{J-infty}\end{equation}
with a uniform constant $C$ using a comparison estimate as in the proof of
Theorem 4.5.1 of \cite{jost}. Applying \eqref{J-infty} separately to $J_1$
and $J_2$, we find that \eqref{J-est} implies the desired estimate.
\end{proof}

We can finally give the proof of Proposition \ref{full-ang-est}, which is
the final geometric lemma required for Theorem \ref{main}.

\begin{proof}[Proof of Proposition \ref{full-ang-est}] In a conic
neighbourhood $U$ of the diagonal this is precisely the content of
Lemma~\ref{lemma:ang}. Outside $U$, we use Proposition~\ref{cor:dist} and
split $d_M(z',z'')$ into the conic distance function $d_{\conic}(z',z'')$
plus a symbolic error $e$.  The smoothness of $e$ on
$\Omega(\eta,r_0)\backslash U$ gives an error term of the form
\begin{equation}
\frac{|v'|^2}{(r')^2} + \frac{|v''|^2}{(r'')^2}.
\label{discr}\end{equation}
As for the conic metric, we need to compare second derivatives of the
distance function along geodesics for $g$, to second derivatives of the
distance function along geodesics for $g_{\conic}$. The geodesics differ by
$O(t^2/r^2)$ as one can see from the geodesic equation
$$
\frac{d^2}{dt^2} z^i = \Gamma^i_{jk} \frac{d}{dt} z^j \frac{d}{dt} z^k
$$
since Christoffel symbols for $g$ differ from Christoffel symbols for $g_{\conic}$ by $O(r^{-2})$ (with respect to a unit length frame). Hence this discrepancy also gives rise to an error of the form \eqref{discr}. Finally the conic calculation gives rise to error terms of the form \eqref{conic-ang-est}. This completes the proof.
\end{proof}


\section{Appendix II: Proof of local smoothing estimates}\label{sec:doi}

In this section we prove the local smoothing estimates \eqref{eq:manifold-doi-comp} -- \eqref{eq:manifold-doi-ang} and Lemma \ref{lemma:superDoi}.  The proof of these
types of estimates in manifolds usually proceeds by the positive commutator
method using pseudo-differential operators which are adapted to the
geometry of geodesic flow.  We shall continue to follow this method,
however instead of using standard pseudo-differential operators, we shall
use the calculus of \emph{scattering pseudo-differential operators} on
$\Mbar$, described by Melrose \cite{melrose}. This calculus is based on the
quantization of certain types of \emph{scattering symbol classes}
$S^{m,l}_{1,\rho}(\overline{M})$, which differ from standard symbols in a
number of ways (in particular, near $\partial M$, every derivative in $z$
of the symbol gains a power of $r = 1/x$), and are defined as follows.

\begin{definition}\label{def:symb} Let $m,l$ be real numbers and $0 \leq \rho < 1$. 
A smooth functions $a: \Tscstar \Mbar \to \C$ is said to be in the \emph{scattering symbol class}
 $S^{m,l}_{1, \rho}(\Mbar)$ provided that one has the usual Kohn-Nirenberg symbol estimates 
$$ | \nabla_z^\alpha \nabla_\zeta^\delta a(z,\zeta) | \leq C_{\alpha,\delta,K} (1 + |\zeta|)^{m-|\delta|}$$
for $z$ ranging in any compact subset $K$ of $M$, and all $\alpha, \beta \geq 0$, as well as the
scattering region estimates
$$
\abs{\partial_x^\alpha \nabla_y^\beta \partial_\nu^\gamma \nabla_\mu^\delta
a(x,y,\nu,\mu)}\leq C_{\alpha,\beta,\gamma,\delta}
\, (1 + |\nu| + |\mu|)^{m-\gamma-\abs\delta}
x^{l-\alpha - \rho(\gamma + \abs\delta)} 
$$ 
whenever we are in the scattering region $(0,\epsilon_0) \times \partial M$, where the dual variables $\nu$, $\mu$
are as in the previous appendix.  Equivalently, using radial and angular derivatives instead of $x$ and $y$ derivatives, we have
$$
\abs{\partial_r^\alpha \nabb^\beta \partial_\nu^\gamma \nabla_\mu^\delta
a(1/r,y,\nu,\mu)}\leq C_{\alpha,\beta,\gamma,\delta}
\, (1 + |\nu| + |\mu|)^{m-\gamma-\abs\delta}
r^{-l-\alpha - |\beta| + \rho(\gamma + \abs\delta)}.
$$ This defines seminorms $\| \cdot \|_{\alpha, \beta,\gamma,\delta}$ on
$S^{m,l}_{1, \rho}(\Mbar)$ in the usual way.  Every symbol $a$ in
$S^{m,l}_{1,\rho}(\Mbar)$ can be quantized to give a \emph{scattering
pseudo-differential operator} $Op(a)$ in the class $\Psisc^{m,l; \rho}(M)$;
the exact means of quantization is not particularly important for our
purposes, but we can for instance use the Kohn-Nirenberg quantization (on
any $n$-dimensional asymptotically conic manifold)
$$
\Op(a) u(z') := (2\pi)^{-n} \int \chi(z',z'')e^{i\langle -\exp^{-1}_{z'}(z''), \zeta \rangle} a(z',\zeta) u(z'')\, dz'' \, d\zeta ,
$$ to define the kernel of $\Op(a)$ near $\partial M$; here, $\chi$ is a
cutoff near the diagonal, chosen so that the inverse of the exponential
map is well defined on $\supp \chi.$
\end{definition}

\begin{example}  If $m, l \in \R$, 
then an operator of the form $(1+H)^{m/2} \langle z\rangle^{-l}$ or
$ \langle z\rangle^{-l} (1+H)^{m/2}$ lies in $\Psisc^{m,l; \rho}$ for any $0 \leq \rho < 1$.  
Heuristically, the quantization of $a(x,y,\mu,\nu)$ can be thought of as 
$a(1/r, y, \frac{1}{i} \partial_r, \frac{1}{i} \nabb)$.
\end{example}

\begin{remark} Scattering symbols $S^{m,l}_{1,\rho}(\Mbar)$ differ from their more standard counterparts $S^{m,0}_{0,0}(\Mbar)$ in that there is a decay of $r^{-l}$ in the symbol, that differentiation in the spatial 
directions gains a power of $r$, and that differentiation in the frequency directions loses a power of $r^\rho$.
The case $\rho = 0$ is the most classical case and is used to prove
\eqref{eq:manifold-doi-comp} -- \eqref{eq:manifold-doi-ang}; however the proof of Lemma~\ref{lemma:superDoi} 
requires the use of more exotic symbols with $0 < \rho < 1$, so that differentiating in the frequency variables $\nu, \mu$ costs us a small power of $r$; this is due to a certain cutoff in $\mu$ which is necessary in our
argument. 
\end{remark}

We now review the calculus for these scattering pseudo-differential operators, obtained\footnote{Strictly speaking, these results were only obtained in \cite{melrose} for $\rho = 0$ but the case for $0 < \rho < 1$ follows from the same argument.} in \cite{melrose}.  If $A$ is an operator in $\Psisc^{m,l; \rho}(\Mbar)$, then its symbol $\sigma(A) \in
S^{m,l}_{1,\rho}(\Mbar)$ is well defined modulo a lower order symbol in $S^{m-1,l+1-\rho}_{1,\rho}(\Mbar)$, and its adjoint $A^*$ is also
in this class with symbol
\begin{equation}\label{symbol-adjoint}
\sigma(A^*) = \overline{\sigma(A)} + O(S^{m-1,l+1-\rho}_{1,\rho}(\Mbar))
\end{equation}
where we use $O(S^{m,l}_{1,\rho}(\Mbar))$ to denote an error in the class $S^{m,l}_{1,\rho}(\Mbar)$.
If $A \in \Psisc^{m,l; \rho}(\Mbar)$ and $B \in \Psisc^{m',l'; \rho}(\Mbar)$ then
\begin{equation}\label{scattering-calculus}
AB \in \Psisc^{m+m',l+l'; \rho}(\Mbar); \quad i[A,B] \in \Psisc^{m+m'-1,l+l'+1 - \rho; \rho}(\Mbar).
\end{equation}
and in fact we have the more precise formulae
\begin{equation}\label{scattering-calculus-symbol}
\begin{split}
\sigma(AB) &= \sigma(A) \sigma(B) + O(S^{m+m'-1,l+l'+1 - \rho}_{1,\rho}(\Mbar))\\
\sigma(i[A,B]) &= \{ \sigma(A), \sigma(B) \} + O(S^{m+m'-2,l+l'+2 - 2\rho}_{1,\rho}(\Mbar))
\end{split}
\end{equation}
where $\{ f, g \}$ is the usual Poisson bracket on the cotangent bundle
$T^* \Mbar$.  Recall that the Poisson bracket $\{f,g\}$ may also be written
as $X_f(g)$ where $X_f$ is the Hamilton vector field associated to $f.$
This identification is crucial in the commutator arguments which follow.

We introduce the weighted Sobolev spaces $H^{m,l}(M)$ as
$$ H^{m,l}(M) := \{ u: \langle z \rangle^l u \in H^m(M) \}.$$
From \eqref{scattering-calculus}, \eqref{scattering-calculus-symbol} and $L^2$-boundedness of $\Psisc^{0,0; \rho}(\Mbar)$ it is easy to verify 
 that operators in $\Psisc^{m,l;\rho}(\Mbar)$ map $H^{m',l'}(M)$ to $H^{m'-m,l'+l}(M)$.

Now let $u(t)$ evolve via the flow \eqref{nls-H}.  Sinced $e^{-itH}$ is unitary and
commutes with $(1+H)^{s/2}$, we have
the trivial estimate
\begin{equation}\label{ls-triv}
 \int_0^1 \| u(t) \|_{H^{s,0}(M)}^2\ dt = \| u(0) \|_{H^s(M)}^2
\end{equation}
for any $s \in \R$.  We now obtain a local smoothing estimate and a
Morawetz estimate for this flow; note that we are now putting the $\ang{z}$
weights into Sobolev norms rather than writing them explicitly:

\begin{lemma}\label{ls-lemma}  For any $s \in \R$ and $\eps > 0$, we have
\begin{equation}\label{ls-scat}
 \int_0^1 \| u(t) \|_{H^{s,-1/2-\eps}(M)}^2\ dt \leq C_{s,\eps} \| u(0) \|_{H^{s-1/2}(M)}^2
\end{equation}
and
\begin{equation}\label{ls-ang}
 \int_0^1 \| \chi \nabb u(t) \|_{H^{s-1,-1/2}(M)}^2\ dt \leq C_s \| u(0) \|_{H^{s-1/2}(M)}^2
\end{equation}
where $\chi$ is a smooth cutoff to the scattering region $M \backslash K_0$.  
\end{lemma}

Observe that this lemma, combined with the above machinery, immediately gives \eqref{eq:manifold-doi-comp} -- \eqref{eq:manifold-doi-ang}.

\begin{proof}  We first observe that it suffices to prove these estimates for a single value of $s$, since
the general case then follows by applying a power of $(1 + H)$ (which commutes with $e^{-itH}$); note
that any error terms arising from commuting $(1+H)$ with, for instance, $\chi \nabb$, can be dealt
with by \eqref{scattering-calculus} and \eqref{ls-triv}.

We begin by proving the weaker estimate \eqref{local-smoothing-manifold}, which in this
language becomes
\begin{equation}\label{ls-manif}
 \int_0^1 \| \varphi u(t) \|_{H^{s,0}(M)}^2\ dt \leq C_{s,\varphi} \| u(0) \|_{H^{s-1/2}(M)}^2
\end{equation}
for any compactly supported non-negative bump function
$\varphi$.  As before it suffices to prove this for a single value of $s$,
say $s = 0$.  This estimate is essentially in \cite{cks}, \cite{doi} and we only give a sketch here.
Using the non-trapping hypothesis of $M$, we may construct a symbol $a(z,\zeta) \in S^{-1,0}_{1,0}(\Mbar)$ such
that $\{ \sigma(H), a \} > 0$ for all $|\zeta|_g \geq 1$ (say), with the uniform bound $\{ \sigma(H), a \} > c > 0$
when $z$ is in the support of $\varphi$.  Note that the Poisson bracket $\{ \sigma(H), a \}$ is nothing more than 
the derivative $\frac{d}{ds} a(z(s),\zeta(s))$
along the geodesic flow \eqref{ham-flow}, which escapes to infinity by hypothesis.  See \cite{cks}, \cite{doi} for
the details of this standard construction.  If we let $A = \Op(a) \in \Psisc^{-1,0;0}(\Mbar)$ be the 
quantization of this operator, then it maps $H^{-1/2}(M)$ to $H^{1/2}(M)$, and thus
$$ |\langle A u(t), u(t) \rangle_M| \leq C \| u(t) \|_{H^{-1/2}(M)}^2 = C \| u(0) \|_{H^{-1/2}(M)}^2$$
for $t=0,1$.  Thus by \eqref{heisenberg} we have
$$ \int_0^1 \langle i[H,A] u(t), u(t) \rangle_M \, dt \leq C \| u(0) \|_{H^{-1/2}(M)}^2.$$
But from the positivity of $\{ \sigma(H), a \}$ we may write
$$ \{ \sigma(H), a \} = c \varphi + |e|^2 + O( S^{-1,0}_{1,0}(\Mbar) )$$
where $e \in S^{0,0}_{1,0}(\Mbar)$ and $c>0$ is a constant.  Quantizing this using \eqref{scattering-calculus-symbol}, 
\eqref{symbol-adjoint} we obtain
$$ i[H,A] = c \varphi + E^* E + O( \Psisc^{-1,0;0}(\Mbar) )$$
for some $E \in S^{0,0}_{1,0}(\Mbar)$, where $O( \Psisc^{m,l;\rho}(\Mbar))$ denotes a scattering pseudo-differential operator in the class $\Psisc^{m,l;\rho}(\Mbar)$.  The term $E^* E$ gives a positive error which can then be discarded,
while the contribution of the $O( \Psisc^{-1,0;0}(\Mbar) )$ error is bounded by \eqref{ls-triv}, and the claim follows.

Now we prove \eqref{ls-ang}.  Here the most convenient value of $s$ is
$s=1$.  Let $a(x,y,\mu,\nu) := \psi^2 \nu$, where $\psi(x,y)$ is a smooth
cutoff to the region $|x| < \epsilon$, and $0 < \epsilon \ll \epsilon_0$ is
a small number to be chosen later; compare with Example \ref{1-particle}.
Let $A$ be the quantization of $a$.  Since $a \in S^{1,0}_{1,0}(\Mbar)$, we
have $A \in \Psisc^{1,0;0}(\Mbar)$, and we use \eqref{heisenberg} as before
to obtain
$$ \int_0^1 \langle i[H,A] u(t), u(t) \rangle_M \leq C \| u(0) \|_{H^{1/2}(M)}^2.$$
From \eqref{ham-flow} we see that
$$ \{ \sigma(H), a \} = \psi^2 (x h^{jk} \mu_j \mu_k) + \frac{\psi^2}{2} x^2 \frac{\partial h^{jk}}{\partial x} \mu_j \mu_k
+ \{ \sigma(H), \psi^2 \} \nu.$$
The second term is dominated by the first in the region $|x| < \epsilon$ if $\epsilon$ is sufficiently small,
while the third term is compactly supported and in $S^{2,0}_{1,0}(\Mbar)$.  We can thus write
$$ \{ \sigma(H), a \} = c \psi^2 x h^{jk} \mu_j \mu_k + |e|^2 + \varphi^2 O(S^{2,0}_{1,0}(\Mbar)) + O(S^{1,0}_{1,0}(\Mbar))$$
for some $e \in S^{1,0}_{1,0}(\Mbar)$ and some compactly supported function $\varphi(z)$.  Quantizing this we obtain
$$ i[H,A] = c (x^{1/2} \nabb_j \psi)^* h^{jk} (x^{1/2} \nabb_k \psi) + E^* E + \varphi O(\Psisc^{2,0;0}(\Mbar)) \varphi + O(\Psisc^{1,0;0}(\Mbar)).$$
The second term is again positive and can be discarded, while the third term can be controlled using
\eqref{ls-manif} and the fourth term by \eqref{ls-triv}.  The claim \eqref{ls-ang} follows (the contribution
of the region $|x| > \epsilon$ being controlled by \eqref{ls-manif}).

Finally we prove \eqref{ls-scat}.  Again we set $s=1$, and now use the symbol $a := -\psi^2 x^{2\eps}
\nu$ where $\psi$ is as
before.  Since $a \in S^{1,0}_{1,0}(\Mbar)$, the quantization $A$ is in $\Psisc^{1,0;0}(\Mbar)$, so by \eqref{heisenberg}
as before we have
$$ \int_0^1 \langle i[H,A] u(t), u(t) \rangle_M \leq C \| u(0) \|_{H^{1/2}(M)}^2.$$
From \eqref{ham-flow} we have
$$ \{ \sigma(H), a \} = 2\eps x^{1+2\eps} \nu^2 \psi^2 - x^{2\eps} \{
\sigma(H), \psi^2 \nu \}.$$ From the proof of \eqref{ls-ang} we see that
$$ \{ \sigma(H), \psi^2 \nu \} = x \psi^2 \mu_j \mu_k O(S^{0,0}_{1,0}(\Mbar)) +
\varphi^2 O(S^{2,0}_{1,0}(\Mbar)) + O(S^{1,0}_{1,0}(\Mbar))$$
for some compactly supported $\varphi$, and thus 
\begin{multline*}
i[H,A] = C \eps (x^{1/2+\eps} \partial_r \psi)^* (x^{1/2+\eps} \partial_r \psi) 
\\ + \nabb^* x^{1/2} \psi O(\Psisc^{0,0;0}(\Mbar)) \psi x^{1/2} \nabb
+ \varphi O(\Psisc^{2,0;0}(\Mbar)) \varphi + O(\Psisc^{1,0;0}(\Mbar)).
\end{multline*}
The second error term is controlled by \eqref{ls-ang}, the third by \eqref{ls-manif}, and the fourth by
\eqref{ls-triv}.  This proves \eqref{ls-scat} in the region $|x| < \epsilon$; the region $|x| > \epsilon$ is
of course controlled by \eqref{ls-manif}.
\end{proof}

To prove Lemma \ref{lemma:superDoi}, it turns out not to be practicable to apply the positive
commutator method directly, mainly because the $z$ integrals in Lemma \ref{lemma:superDoi}
are not themselves positive.  Instead, we use the positive commutator method to first
estimate an auxiliary positive quantity, which we then use to control the expressions in
Lemma \ref{lemma:superDoi}.  More specifically, we prove the following variant of \eqref{ls-ang} which
only uses `half' an angular derivative instead of a full angular derivative.

\begin{lemma}[Half-angular Morawetz estimate]\label{bpm}  
Let $0 < \rho < 1$, and let $\phi: \R \to \R$ be a smooth non-decreasing function such that $\phi = 0$ on
$(-\infty,1]$ and $\phi=1$ on $[2,+\infty)$.  Let $b(x,y,\mu,\nu)$ be the $S^{1,1/2}_{1,\rho}(\Mbar)$
symbol
\begin{equation}
b := \frac1{2} \phi(\frac{\epsilon_0}{x})
\phi(|\mu|^2 + |\nu|^2) 
\phi\big(\frac{\abs{\mu}}{|\nu| x^\rho}\big) x^{1/2} 
\Big( \phi( \frac{|\nu|}{|\mu|} )  |\mu|^{1/2} |\nu|^{1/2} + |\mu|  \Big)
\label{b-defn}\end{equation}
where $|\mu|^2 = |\mu|^2_{h(x)} := h^{jk}(x,y) \mu_j \mu_k$.  Let $B \in \Psisc^{1,1/2;\rho}(\Mbar)$ be a
quantization of this operator.  Then for any $H^{s+1/2}$ solution to \eqref{nls-H} and any $s \in \R$, we have
\begin{equation}
\label{commutatorestimate}
\int_0^1 \| B u(t) \|_{H^s(M)}^2\ dt \leq C
\norm{u(0)}_{H^{s+1/2}(M)}^2.
\end{equation}
Here the constants depend on $\rho$ and $s$.
\end{lemma}

\begin{remark} Ignoring all the cutoffs, the important part of $b$ is the $|\mu|^{1/2} |\nu|^{1/2}$ term.  Heuristically, $B$ can be thought of as the operator $ r^{-1/2} |\nabb|^{1/2} |\partial_r|^{1/2}$,
with the various cutoffs involving $\phi$ needed to avoid a singularity arising from any degeneracy of
$\nabb$ or $\nabla$, localizing physical space to the scattering region $r > 1/\epsilon_0$, and frequency space
to the region $|\nabb| \geq |\nabla|  r^{-\rho} \geq r^{-\rho}$.  This cutoff is not dangerous for us as
in the region where $|\nabb| \leq |\nabla| r^{-\rho}$, one can effectively deduce Lemma \ref{lemma:superDoi}
from \eqref{ls-scat} (this is why we require $\rho > 0$).  The reader may wish to 
ignore the presence of the $\phi$ cutoffs for a first reading.  The estimate \eqref{commutatorestimate} can then 
be thought of, very non-rigorously, as
$$ \int_0^1 \| \chi |\nabla|^{1/2} |\nabb|^{1/2} u(t) \|_{H^{s,-1/2}(M)}^2\ dt
\leq C \| u(0) \|_{H^{s+1/2}(M)}^2$$
for some appropriate cutoff $\chi$; compare this with \eqref{ls-ang}.  The crucial point 
here is that we do not lose an epsilon of decay as we would from \eqref{ls-scat}.
\end{remark}

\begin{proof}
As in the proof of Lemma \ref{ls-lemma} it suffices to verify this when $s=0$.  We again
use the positive commutator method.  Morally speaking, the commutant to use 
is $|\nabb| \sgn(\frac{1}{i} \partial_r)$; however this operator is too singular, and so we must use a modified version of this commutant.

Let $a$ denote the function \begin{equation*}
a(x, y, \mu, \nu) := \phi^2(\epsilon_0/x) \phi^2(\nu^2 + |\mu|^2)  \bigg( \phi^2\big(\frac{|\mu|}{|\nu| x^\rho}\big)  \phi^2(\frac{|\nu|}{|\mu|}) (-|\mu| \sgn \nu) + C \nu \bigg),
\end{equation*}
where $C$ is a large constant to be chosen later. 
Notice that $a$ vanishes in the near region $K_0$.
This can easily be verified to be a scattering symbol of order $S^{1,0}_{1,\rho}(\Mbar)$, because of all the cutoffs.  Let $A \in \Psisc^{1,0;\rho}(\Mbar)$ be the quantization of $a$.  
Applying \eqref{heisenberg} as in Lemma \ref{ls-lemma} we have
$$ \int_0^1 \langle i[H,A] u(t), u(t) \rangle_M \leq  C \| u(0) \|_{H^{1/2}(M)}^2.$$
We shall now factorize $i[H,A]$ in the form
\begin{equation}
 i[H,A] = B^* B  + \sum_{\text{finite}} C_i^* C_i + O( \Psisc^{2,1+\delta;\rho}(\Mbar)) + O( \Psisc^{1,0;\rho}(\Mbar) )
 \label{poscomm}\end{equation}
where $C_i \in \Psisc^{1,1/2;\rho}(\Mbar)$ and   $\delta > 0$; this will prove the claim since the error terms can be treated by \eqref{ls-scat}, \eqref{ls-triv}.
From \eqref{scattering-calculus-symbol} it suffices to obtain a representation
\begin{equation}
 \{ \sigma(H), a \} = |b|^2 + \sum_i |c_i|^2 +  O(S^{2,1+\delta}_{1,\rho}(\Mbar)) + O( S^{1,0}_{1,\rho}(\Mbar)).
\label{repr}\end{equation}

Observe that $\{ \sigma(H), a \}$ vanishes when $x > \epsilon_0$, and in the region $\epsilon_0 \geq x \geq \epsilon_0/4$
it is a symbol of order $S^{2,1-\rho}_{1,\rho}(\Mbar)$ as mentioned earlier, and hence in this region is also
a symbol of order $S^{2,l}_{1,\rho}(\Mbar)$ for any $l$.  Thus this factorization is easy to obtain in the near
region, and we can focus on the scattering region $x < \epsilon_0/2$.  In particular we can now ignore
the cutoff $\phi(\epsilon_0/x)$.  A similar argument allows us to work in the region where
$|\mu|^2 + |\nu|^2 > 2$ (as $\{ \sigma(H),a \}$ is certainly in $S^{1,0}_{1,\rho}(\Mbar)$ when $|\mu|^2 + |\nu|^2 < 4$), allowing us to ignore the $\phi(|\mu|^2 + |\nu|^2)$ cutoff.  We thus need to consider 
$$
\{ \sigma(H),  \phi^2\big(\frac{|\mu|}{|\nu| x^\rho}\big)  \phi^2(\frac{|\nu|}{|\mu|}) (-|\mu| \sgn \nu) + C \nu \}
$$
in the region $x < \epsilon_0/2$, $|\mu|^2 + |\nu|^2 > 2$. 

Since the two $\phi^2$ terms are each $1$ on the support on the derivative of the other, we may express this as 
\begin{equation}\begin{gathered}
\phi^2\big(\frac{|\mu|}{|\nu| x^\rho}\big)  \phi^2(\frac{|\nu|}{|\mu|})
\{ \sigma(H),   -|\mu| \sgn \nu  \}  + C \{ \sigma(H),   \nu \} \\
 -|\mu| \sgn \nu \, \{ \sigma(H),    \phi^2(\frac{|\nu|}{|\mu|})  \} 
 -|\mu| \sgn \nu \, \{ \sigma(H),  \phi^2\big(\frac{|\mu|}{|\nu| x^\rho}\big) \}  .\end{gathered}\end{equation}
We have $\{ \sigma(H),   -|\mu| \sgn \nu  \} = x \big(|\mu| |\nu| + O(S^{2,1}_{1,\rho}(\Mbar)) \big)$ and $\{ \sigma(H),   \nu \} = x \big(|\mu|^2 + O(S^{2,1}_{1,\rho}(\Mbar)) \big)$. Hence for $C > 1$ the sum of these terms is larger than $|b|^2$, where $b$ is as in \eqref{b-defn}. To analyse the third term we compute
$$
\{ \sigma(H),    \phi^2(\frac{|\nu|}{|\mu|})  \} = 2x \phi \phi'(\frac{|\nu|}{|\mu|}) \big( |\mu| + \frac{\nu^2}{|\mu|} \big) + O(S^{2,2}_{1,\rho}(\Mbar)).
$$
This term is dominated by the second term, modulo acceptable errors, for sufficiently large $C$. As for the fourth term, we compute 
\begin{equation}\begin{gathered}
|\mu| \{ \sigma(H), \phi^2(\frac{|\mu|}{x^\rho |\nu|}) \} = 
-2 x  |\mu|  \phi\phi'(\frac{|\mu|}{ x^\rho |\nu|}) 
\Big( \frac{(1 + \rho) |\mu|}{x^\rho} + \frac{|\mu|^3}{x^\rho \nu^2 } + O((S^{1,1}_{1,\rho}(\Mbar)) \Big).
\end{gathered}\label{cutoff2}\end{equation}
These terms are all supported where $|\mu|$ is comparable to $x^\rho \nu$. Since there is an overall factor of $|\mu|^2$ on the right hand side of \eqref{cutoff2}, we gain a factor of $x^{2\rho}$, showing that  \eqref{cutoff2} is a symbol of order $(2, 1 + \rho)$ which is an error term in \eqref{repr}. 

Assume that the function $\phi$ is such that terms such as $\sqrt{\phi}$ and $\sqrt{\phi'}$, etc, are smooth; for example, we make take
$$
\phi(t) = 1 - \phi(3 - t) \text{ and } \phi(t) = e^{-(t-1)^{-1}} \text{ for } t \in [1, 1.2].
$$
Then $\{ \sigma(H), a \}$ can be written, modulo acceptable errors, as $|b|^2$ plus a sum of squares $\sum |c_i|^2$. This completes the proof of \eqref{poscomm}, which establishes the lemma. 
\end{proof}

\begin{proof}[Proof of Lemma \ref{lemma:superDoi}]
We first prove \eqref{super-1}.  Observe if the derivative $\nabla_k$ was replaced by an angular derivative
then this claim would follow from \eqref{eq:manifold-doi-ang} or \eqref{ls-ang}, so it suffices to show that
$$
\int_0^1 \sup_{w \in W} \bigg| \int_{M} \frac{a^{j}_w(z) \nabb_j u(t,z)
\overline{\partial_r u(t,z) }}{r} \, dg(z) \bigg| dt \leq C \| u(0)
\|_{H^{1/2}(M)}^2
$$
where $a^j_w$ satisfies similar symbol estimates to $a^{jk}_w$.
Choose a  $0 < \rho < 1$ (for instance, we may take $\rho := 1/2$).
Observe that we may write
$$ \int_{M} \frac{a^{j}_w(z) \nabb_j u(t,z) 
\overline{\partial_r u(t,z) }}{r} \, dg(z) = \int A_w u(t,z)
\overline{u(t,z)}\, dg(z)$$
where $A_w \in \Psisc^{2,1;\rho}(\Mbar)$ has symbol
$$ \sigma(A_w)(x,y,\nu,\mu) = C a^{j}_w(x,y) x \mu_j \nu + O( S^{1,0}_{1,\rho} )$$
for some constant $C$ 
(indeed we may improve the error term substantially). It is here that we require the functions $a^j_w$ to be symbols; otherwise $A_w$ would not be a pseudodifferential operator.  Observe 
that in the region $|\mu|^2 + |\nu|^2 \leq 4$, the main term is in $S^{1,0}_{1,\rho}$.  From the support
hypothesis of $a^j$ we thus have
$$ \sigma(A_w)(x,y,\nu,\mu) = C \phi(\epsilon_0/x)
\phi(|\mu|^2 + |\nu|^2) a^{j}_w(x,y) x \mu_j \nu + O( S^{1,0}_{1,\rho} ).$$
Recalling the operator $B \in \Psisc^{1,1/2;\rho}(\Mbar)$ with symbol $b \in S^{1,1/2}_{1,\rho}(\Mbar)$ from
the previous lemma, we can thus factorize
\begin{multline*} \sigma(A_w)(x,y,\nu,\mu) = O(S^{0,0}_{1,\rho}(\Mbar)) |b|^2
\\+ C \phi(\epsilon_0/x) \phi(|\mu|^2 + |\nu|^2) (1 - \phi\big(\frac{|\mu|}{|\nu| x^\rho}\big))
a^{j}_w(x,y) x \mu_j \nu + O( S^{1,0}_{1,\rho}(\Mbar)).
\end{multline*}
The second term is supported on the region where $|\mu| \leq 2 |\nu| x^\rho$ and can thus
be easily seen to lie in $S^{2,1+\rho}_{1,\rho}(\Mbar)$.  Thus by \eqref{scattering-calculus-symbol} we have
$$ A_w = B^* O(\Psisc^{0,0;\rho}(\Mbar)) B +
O( \Psisc^{2,1+\rho;\rho} ) + O( \Psisc^{1,0;\rho} );$$
since $\Psisc^{m,l;\rho}$ maps $H^{m/2,-l/2}(M)$ to $H^{-m/2,l/2}(M)$,
$$ \big| \int A_w u(t,z)\overline{u(t,z)} \, dg(z) \big| \leq C \| B u(t) \|_{L^2(M)}
+ C \| u(t) \|_{H^{1,-\frac{1+\rho}{2}}(M)} + C \| u(t) \|_{H^{\frac1{2}}(M)}$$
uniformly in $w$.  The claim \eqref{super-1} then follows from Lemma
\ref{bpm}, \eqref{ls-scat} and \eqref{ls-triv}.

Now we prove \eqref{super-2}.  As before we have
$$ \int_{M} \frac{a^{j}_w(z) \nabb_j u(t,z) 
u(t,z)}{r} \, dg(z) = \int \tilde A_w u(t,z)\overline{u(t,z)}dg(z)$$
where $\tilde A_w \in \Psisc^{1,1;\rho}(\Mbar)$ has symbol
$$ \sigma(\tilde A_w)(x,y,\nu,\mu) = C a^{j}_w(x,y) x \mu_j + O( S^{0,0}_{1,\rho} ),$$
which is similar to the previous but without the factor of $\nu$.
As before, we localize the main term to the region $|\mu|^2 + |\nu|^2 \geq 1$.  We
can then smoothly split the main term into the region where $|\nu| \geq \frac{1}{2} |\mu|$
and where $|\mu| \geq |\nu|$.  In the first region, we have $(1 + \sigma(H))^{1/2}$ comparable
to $|\nu|$, and so by repeating the previous argument we eventually obtain
$$ \tilde A_w = B^* O(\Psisc^{-1,0;\rho}(\Mbar)) B +
O( \Psisc^{1,1+\rho;\rho} ) + O( \Psisc^{0,0;\rho} ),$$
and one then argues as before.
In the second region, we have $(1 + \sigma(H))^{1/2}$ comparable to $|\mu|$, so we can write directly
$$ \sigma(\tilde A_w)(x,y,\nu,\mu) = O( S^{-1, 0}_{1,\rho} ) x |\mu|^2 + O( S^{0,0}_{1,\rho} ),$$
and hence we have a decomposition
$$ \tilde A_w = \nabb O(\Psisc^{-1,0;\rho}(\Mbar)) \nabb + O( \Psisc^{0,0;\rho} )$$
uniformly in $w$.  The claim now follows using \eqref{ls-ang} instead of Lemma \ref{bpm}.
\end{proof}


\end{document}